\documentclass[1p,number,times]{elsarticle}
\usepackage[pagewise]{lineno}
\usepackage{amsmath, amsfonts,amsthm, amssymb}
\usepackage{mathrsfs}
\usepackage[utf8]{inputenc}
\usepackage{xcolor}
\usepackage{graphicx}
\usepackage{mathtools}
\usepackage{cleveref}
\usepackage[all]{xy}
\usepackage[normalem]{ulem}
\usepackage{algorithm}
\usepackage{algpseudocode}
\usepackage[defaultlines=4]{nowidow}
\usepackage{nicematrix}
\usepackage{wasysym}
\usepackage{marvosym}
\usepackage{microtype}
\usepackage[commandnameprefix=always]{changes}
\definechangesauthor[name=Patrice, color=blue]{pom}
\newbox\myboxa
\newbox\myboxb
\def\clap#1{\hbox to 0pt{\hss#1\hss}}

\makeatletter
\def\word#1{\mbox{\@tfor \i:= #1\do{%
\if.\i\relax\makebox[25pt]{$\dots$}\else\makebox[10pt]{$\i$}\fi
}}\ }

\theoremstyle{remark}
\newtheorem{claim}{$\rhd$ Claim}
\newcommand{\cqed}{\ensuremath{\lhd}}
\newenvironment{claimproof}{\par
	\pushQED{\cqed}%
	\normalfont \topsep6\p@\@plus6\p@\relax
	\trivlist
	\item\relax
	{\itshape
		Proof of the claim\@addpunct{.}}\hspace\labelsep\ignorespaces
}{%
	\hfill\popQED\endtrivlist\@endpefalse
}
\makeatother
\newtheorem{theorem}{Theorem}
\newtheorem{corollary}[theorem]{Corollary}
\newtheorem{lemma}[theorem]{Lemma}

\newtheorem{remark}[theorem]{Remark}
\newtheorem{definition}[theorem]{Definition}

\newtheorem{fact}[theorem]{Fact}

\newtheorem{example}[theorem]{Example}
\crefname{lemma}{Lemma}{Lemmas}

\newcommand{\prop}[1]{\ensuremath{\mathtt{#1}}}

\DeclareMathOperator{\Sym}{{\rm Sym}}
\newcommand{\del}{\mathbin{\setminus}}
\newcommand{\cont}{\mathbin{/}}

\newcounter{dummyc}
\renewcommand{\thefootnote}{\Cross}
\journal{European Journal of Combinatorics (special issue)}
\begin{document}
	\begin{frontmatter}
\title{A few words about maps\tnoteref{ERC}}
\tnotetext[ERC]{\ERCagreement}
\author{Robert Cori}\address{Labri, Université Bordeaux 1}\ead{robert.cori@labri.fr}
\author{Yiting Jiang}\address{Universit\'e de Paris, CNRS, IRIF, F-75006, Paris, France and Department of Mathematics, Zhejiang Normal University, China}\ead{yjiang@irif.fr}
\author{Patrice Ossona de Mendez}\address{Centre d'Analyse et de Math\'ematiques Sociales (CNRS, UMR 8557), Paris, France and Computer Science Institute of Charles University, Praha, Czech Republic}\ead{pom@ehess.fr}
\author{Pierre Rosenstiehl\textsuperscript{\Cross}}
	\footnotetext{This project started at the occasion of the talk presented in 2009 by Pierre Rosenstiehl in Bordeaux. Pierre left before this paper was eventually written, but it is fair to include him as a coauthor.}

\newcommand{\ERCagreement}{This paper is part of a project that has received funding from the European Research Council (ERC) under the European Union's Horizon 2020 research and innovation programme (grant agreement No 810115 -- {\sc Dynasnet}) and from the french ANR project HOSIGRA (ANR-17-CE40-0022).\\
	\includegraphics[width=.25\textwidth]{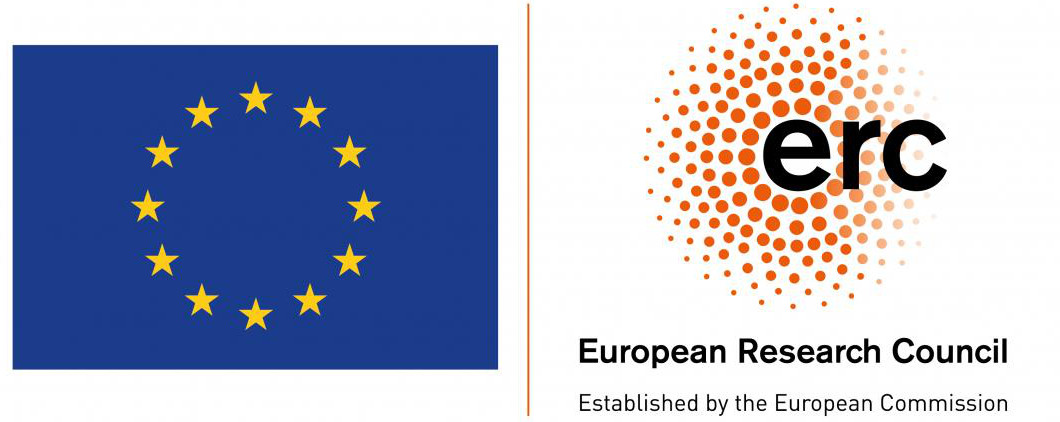}}

	\begin{abstract}
In this paper, we survey some properties, encoding, and bijections involving combinatorial maps, double occurrence words, and chord diagrams. We particularly study quasi-trees from a purely combinatorial point of view and derive a topological representation of maps with a given spanning quasi-tree using two fundamental polygons, which extends the representation of planar maps based on the equivalence with bipartite circle graphs. Then, we focus on Depth-First Search trees and their connection with a poset we define on the spanning quasi-trees of a map. We apply the bijections obtained in the first section to the problem of enumerating loopless rooted maps. Finally, we return to the planar case and discuss a decomposition of planar rooted loopless maps and its consequences on  planar rooted loopless map enumeration.
	\end{abstract}
\renewcommand*{\thefootnote}{\arabic{footnote}}
\end{frontmatter}

\section{Introduction}
Combinatorial maps is a natural bridge linking combinatorics to topological graph theory. At the heart of this connection is the property that a cellular embedding of a graph on a surface can (basically) be described by means of the cyclic order of  incident edges at vertex. This property, already noticed by Heffter \cite{Heffter}
at the end of the XIXth century,  led to a combinatorial study of graph embeddings, initiated  by Edmonds \cite{edmonds} and Youngs \cite{Youngs1}, who independently gave a precise description of this correspondence, as well as a derived computation of the faces and genus of the embedding. This framework was then popularized by White's book ``Graphs, Group, and Surfaces'' \cite{ggs}. Further generalization have been proposed,  with the notions of hypermaps  \cite{Cori2,jacques1968genre},  of dessins d'enfants \cite{schneps1994grothendieck}, and  of ribbon graphs \cite{reshetikhin1990ribbon}. Also, instead of a description of combinatorial maps by means of a permutation and an involution, Tutte proposed a description based on three involutions.

From these combinatorial descriptions of graphs on surfaces,  some coding schemes of maps by words were derived, which allowed to shed a new light in \cite{Cori2} on Tutte's enumeration formulas of maps \cite{tutte1963census}, or a bijection between (rooted) hypermaps and indecomposable permutations (also called connected or irreducible) \cite{newperm} (see also \cite{cori2009indecomposable}).  
The connection between maps and words also impacted several over fields, from graph drawing \cite{PRCpolhor} to matroid theory \cite{Taxi_occu,Taxi_circle}, and was central to the characterization of Gauss codes \cite{GaussPR,Taxi_TheGauss}.

In this paper, we survey some aspects of combinatorial maps, and how they are intrinsically linked with chord diagrams, double occurrence words, Euler tours, $\Delta$-matroids, polynomials, etc.,  with an emphasis on words. In particular, we explicit few bijections involving rooted maps and words.
As an application, we derive some enumerative properties of  loopless rooted maps. Our approach will be purely combinatorial, though it relates to topological properties of graph embeddings on orientable surfaces. The reader interested in the topological aspects is referred to the monograph
\cite{lando2004graphs} (See also \cite{moharthom}).

\section{Combinatorial maps and hypermaps}
Recall that  the group of all permutations of a set $B$ is the \emph{symmetric group} $\Sym(B)$ of $B$.
The composition of permutations is denoted multiplicatively: $\sigma\tau=\sigma\circ\tau$.
A permutation $\tau\in\Sym(B)$ is an \emph{involution} if $\tau^2$ is the identity permutation; it is \emph{fixed-point free} if $\tau(x)\neq x$ for every $x\in B$. We denote by   $\langle\tau_1,\dots,\tau_k\rangle$
the subgroup of $\Sym(B)$ generated by the permutations $\tau_1,\dots,\tau_k$, which is the smallest subgroup of $\Sym(B)$ which contains $\tau_1,\dots,\tau_k$. The \emph{orbit} of an element $b\in B$ in a subgroup $\Gamma$ of $\Sym(B)$ is the set $\Gamma\cdot b=\{\gamma(b): \gamma\in\Gamma\}$.
A subgroup $\Gamma$ of $\Sym(B)$ \emph{acts transitively} on $B$ if, for every $b,b'\in B$ there exists some $\tau\in\Gamma$ such that $\tau(b)=b'$. In other words, $\Gamma$ acts transitively on $B$ if 
$\Gamma$ has a single orbit.
The \emph{cycles} of a permutation $\gamma\in\Sym(B)$ are the orbits of $\langle\gamma\rangle$, that is, the sets of the form $\{\gamma^i(b): i\in\mathbb Z\}$. The set of all the cycles of $\gamma$ is denoted by $Z(\gamma)$. 

For $\mu\in{\rm Sym}(B)$ and $B'\subseteq B$ we define the \emph{restriction}  of $\sigma$ to $B'$
as the permutation $\sigma_{|B'}\in{\rm Sym}(B')$ defined as follows: for every $b\in B'$,
$\mu_{|B'}(b)=\mu^k(b)$, where $k$ is the minimum positive integer such that $\mu^k(b)\in B'$. We further define the \emph{cutting-out} of $\mu$ on $B'$ as the permutation $\mu_{B'}\in{\rm Sym}(B)$, where $\mu_{B'}(b)$ is equal to $\mu_{|B'}(b)$ if $b\in B'$, and to $b$ if $b\notin B'$. (Intuitively, $\mu_{B'}$ is $\mu_{|B'}$ on $B'$ and the identity mapping on $B\setminus B'$.)

The next easy lemma will be useful.
\begin{lemma}
	\label{lem:prodrest}
		Let $\sigma,\mu\in{\rm Sym}(B)$ and let $B'\subseteq B$.
		Then, 
\[\sigma_{B'}\mu_{B'}=(\sigma\mu_{B'})_{B'}\quad\text{and thus}\quad \sigma_{|B'}\mu_{|B'}=(\sigma\mu_{B'})_{|B'}.\]
\end{lemma}
\begin{proof}
	If $b\notin B'$, then $b$ is a fixed point of both $\sigma_{B'}\mu_{B'}$ and $(\sigma\mu_{B'})_{B'}$.
	So, assume $b\in B'$. Let $k$ be the minimum positive integer with $(\sigma\mu_{B'})^k(b)\in B'$.
	We prove by induction on $i$ that for every $1\leq i\leq k$ we have
	$(\sigma\mu_{B'})^i(b)=\sigma^i\mu_{B'}(b)$. The base case, $i=1$, is straightforward.
	Assume that the property holds for all $1\leq j<i$ and $1<i\leq k$.
	Then, 
	\begin{align*}
	(\sigma\mu_{B'})^i(b)&=\sigma\mu_{B'}(\sigma\mu_{B'})^{i-1}(b)\\
	&=\sigma(\sigma\mu_{B'})^{i-1}(b)&\text{(as $(\sigma\mu_{B'})^{i-1}(b)\notin B'$)}\\
	&=\sigma(\sigma^{i-1}\mu_{B'})(b)&\text{(by induction hypothesis)}\\
	&=\sigma^i\mu_{B'}(b)
	\end{align*} 
Thus, it follows that $k$ is also the minimum positive integer with $\sigma^k\mu_{B'}(b)\in B'$.
Hence, $(\sigma\mu_{B'})_{B'}(b)=(\sigma\mu_{B'})^k(b)=\sigma^k\mu_{B'}(b)=\sigma_{B'}\mu_{B'}(b)$.
\end{proof}


A \emph{(combinatorial) map} is a 
triple $\mathcal M=(B,\sigma,\alpha)$, where  $B$ is the set of the \emph{flags} of $\mathcal M$,  $\sigma,\alpha\in{\rm Sym}(B)$,  $\alpha$ is a fixed-point free involution, and $\langle\sigma,\alpha\rangle$ acts transitively on $B$. 
When we relax the condition that  $\langle\sigma,\alpha\rangle$ acts transitively on $B$, we say that $\mathcal M=(B,\sigma,\alpha)$ is a \emph{general map}, and we refer to the orbits of $\langle \sigma,\alpha\rangle$ as the \emph{components} of $\mathcal M$.
The \emph{edges} of a map $\mathcal M$ are the cycles of $\alpha$, and its \emph{vertices}  are the cycles of $\sigma$.
For a flag $b\in B$ it will be convenient to define
 \[\underline{b}=\langle\alpha\rangle\cdot b=\{b,\alpha(b)\},\]
which is  the edge of $\mathcal M$ that contains $b$.

Two maps $\mathcal M=(B,\sigma,\alpha)$ and $\mathcal M'=(B',\sigma',\alpha')$ are \emph{isomorphic} if there exists a bijection $f:B\rightarrow B'$ (referred to as a \emph{(map) isomorphism} of $\mathcal M$ and $\mathcal M'$) with
$\sigma'=f\circ\sigma\circ f^{-1}$ and $\alpha'=f\circ\alpha\circ f^{-1}$. 

The \emph{underlying graph} of a map $\mathcal M$ is a graph  $G(\mathcal M)$ with vertex set $Z(\sigma)$ and edge set $Z(\alpha)$. In this graph,  an edge $e\in Z(\alpha)$ is incident to a vertex $v\in Z(\sigma)$ if $e\cap v\neq\emptyset$. If $e\subseteq v$ (i.e. if $e$ and $v$ intersects on two flags), then $e$ is a \emph{loop} attached to $v$. 
It will be convenient to identify $\mathcal M$ with an embedding of  $G(\mathcal M)$. When we  speak about
cycles,  trees, and cuts of $\mathcal M$, we mean the   cycles,  trees, and cuts of $G(\mathcal M)$, as subsets of edges.

It is well known that a  map $\mathcal M$ defines a cellular embedding of $G(\mathcal M)$ on an orientable surface, whose genus $g$ is given by Euler's formula
\begin{equation}
	2-2g=|Z(\sigma\alpha)|-|Z(\alpha)|+|Z(\sigma)|.
\end{equation}

\begin{figure}[ht] 
	\centering
	\includegraphics[width=.3\linewidth]{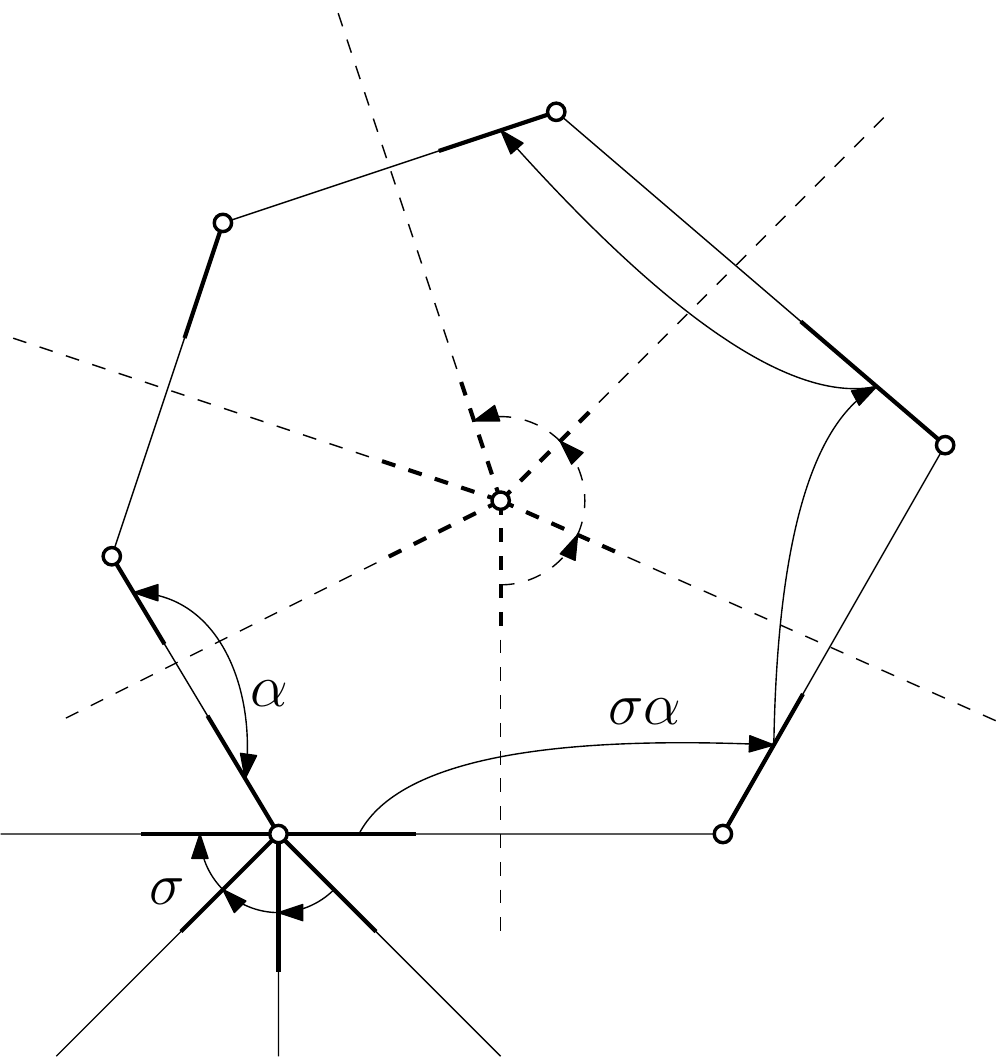}
	\caption{Map duality. From a geometric point of view, the orientation of the surface is  reversed when we consider the dual map.}
	\label{fig:dual}
\end{figure}

In this paper, by the \emph{dual} of a general map $\mathcal M=(B,\sigma,\alpha)$ we mean the general map 
\[\mathcal M^\ast=(B,\sigma\alpha,\alpha).\]
Note that $\mathcal M^\ast$ has the same components as $\mathcal M$; in particular, the dual of a map is a map.
Remark  that this duality  differs from the geometric dual by the orientation of the dual surface (see \Cref{fig:dual}).

A \emph{rooted map} is a pair $\mathcal M_\bullet=(\mathcal M,b_\bullet)$, where $b_\bullet$ is a flag of $\mathcal M$. Two rooted combinatorial maps $(\mathcal M,b_\bullet)$ and $(\mathcal M',b_\bullet')$ are \emph{isomorphic} if there exists a map isomorphism $f$ of $\mathcal M$ and $\mathcal M'$ with  $b_\bullet'=f(b_\bullet)$.

For an introduction to combinatorial maps (and hypermaps) we refer the interested reader to~\cite{Cori}.

\section{Planar maps}
\begin{figure}[ht] 
	\centering
	\includegraphics[width=\linewidth]{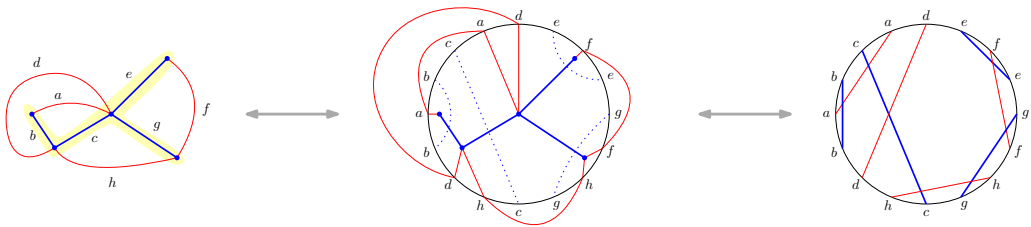}
	\caption{From a spanning tree $T$ of a planar graph $G$ to the bipartite circle graph $\Lambda_T(G)$ , and back.}
	\label{fig:bipcirc}
\end{figure}

Before we consider the general setting of maps, we take time to comment on the planar case.
It follows immediately from Euler's formula that a map $\mathcal M$ is planar if and only if  the complement of a spanning tree of $\mathcal M$ is  a spanning tree of the dual map $\mathcal M^\ast$.
Fixing a spanning tree $T$ of a planar map $\mathcal M$ allows defining several bijections.
A first one is obtained by considering a tubular neighborhood of $T$. Traversing this tubular neighborhood, either we follow some edge of $T$, or we cross some edge of the complement of $T$. This way, the traversal defines a circular sequence, in which every edge appears twice. This circular sequence, in turn, naturally defines a bipartite chord diagram (hence a bipartite circle graph). Conversely, every proper coloration of a bipartite circle graph uniquely determines a planar map and a spanning tree of it (see \Cref{fig:bipcirc}).
This property is at the heart of Fraysseix's characterization of circle graphs \cite{Taxi_occu,Taxi_circle} and Rosenstiehl's characterization of planarity by the algebraic diagonal \cite{rosenstiehl1976characterization}.

For an edge $f\notin T$ we denote by $\gamma_T(f)$ the \emph{fundamental cycle} of $f$, that is the edge set of the unique cycle in $T\cup\{f\}$; 
For an edge $e\in T$ we denote by $\omega_T(e)$ the \emph{fundamental cocycle} of $e$, that is the unique inclusion minimal cut of $G$ intersecting $T$ only at $e$.
It is easily proved that for an edge $e\in T$ and an edge $f\notin T$ we have
$e\in\gamma_T(f)\iff f\in\omega_T(e)$.  The \emph{fundamental interlacement graph} of the tree $T$ in the  graph $G$
if the bipartite graph $\Lambda_T(G)$ with vertex set  $T\cup (E(G)\setminus T)$, where $e\in T$ is adjacent to $f\in E(G)\setminus T$ if $e\in\gamma_T(f)$.  The above described construction emphasizes that the fundamental interlacement graph $\Lambda_T(G)$ of a spanning tree $T$ in a planar graph $G$ is a bipartite circle graph.

\begin{figure}[th]
	\centering
	\includegraphics[width=0.5\linewidth]{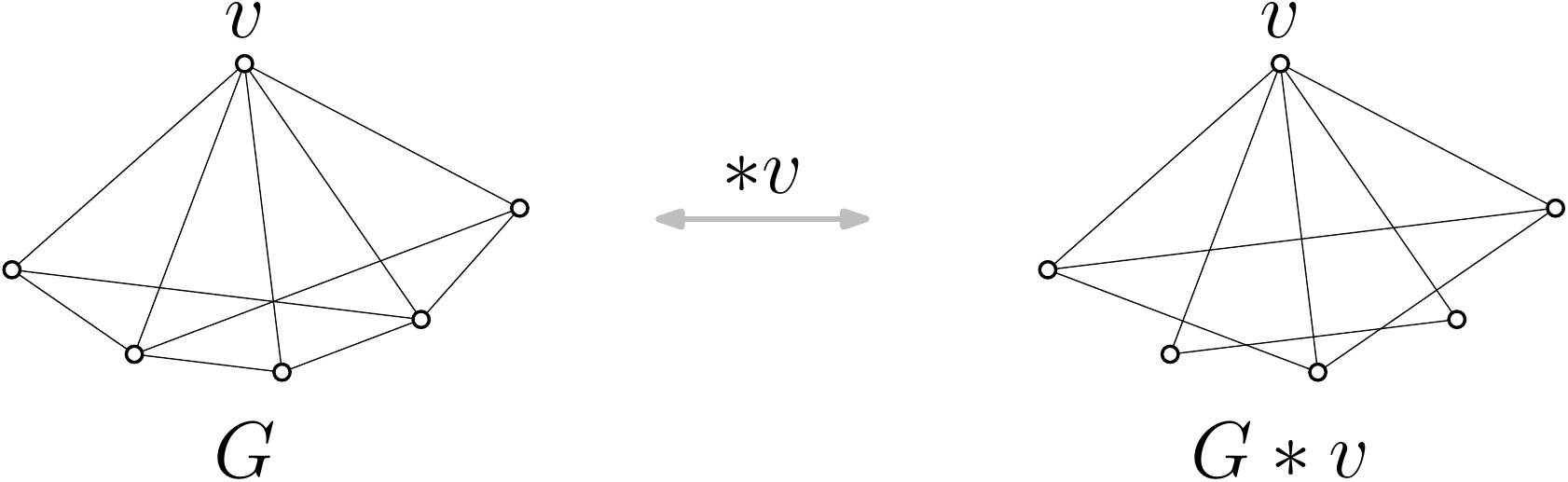}
	\caption{Local complementation of a graph $G$ at a vertex $v$}
	\label{fig:lc}
\end{figure}

The {\em local complementation} of a  graph $G$ at a vertex $v$ is the  graph $G\ast v$  obtained by replacing the subgraph induced by $G$ on the neighbors of $v$ by its complement (see \Cref{fig:lc}). 

This transformation was introduced by Kotzig \cite{kotzig1977quelques} in relations with $\kappa$-transformations.
The \emph{pivoting} (or \emph{switching}) of an edge $uv$ in a graph $G$ results in the graph $G\wedge uv=G\ast u\ast v\ast u=G\ast v\ast u\ast v$.

\begin{figure}[ht] 
	\centering
	\includegraphics[width=\linewidth]{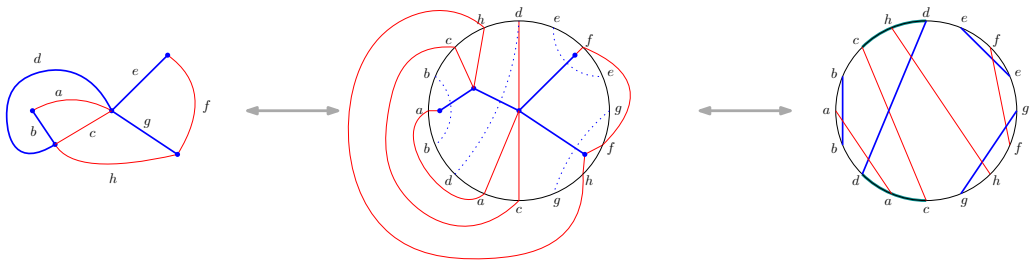}
	\caption{From a planar graph to a bipartite circle graph, and back\dots with another tree obtained by replacing $c\in T$ by $d\notin T$. The corresponding chord diagram is obtained by pivoting $cd$. }
	\label{fig:bipcirc2}
\end{figure}

The chord diagrams and circle graphs obtained when one consider a different spanning tree are nicely related. Indeed, if a spanning tree $T'$ is obtained from a spanning tree $T$ by removing a tree edge $e$ and replacing it with a non-tree edge $f$, then $f\in\omega(e)$ (that is: $ef$ is an edge of $\Lambda_T(G)$), and we have $\Lambda_{T'}(G)=\Lambda_{T}(G)\wedge ef$ (that is: the circle graph  obtained from $\Lambda_T(G)$ by pivoting the edge $ef$; see \Cref{fig:bipcirc2}). 

\begin{figure}[ht]
	\centering
	\includegraphics[width=0.7\textwidth]{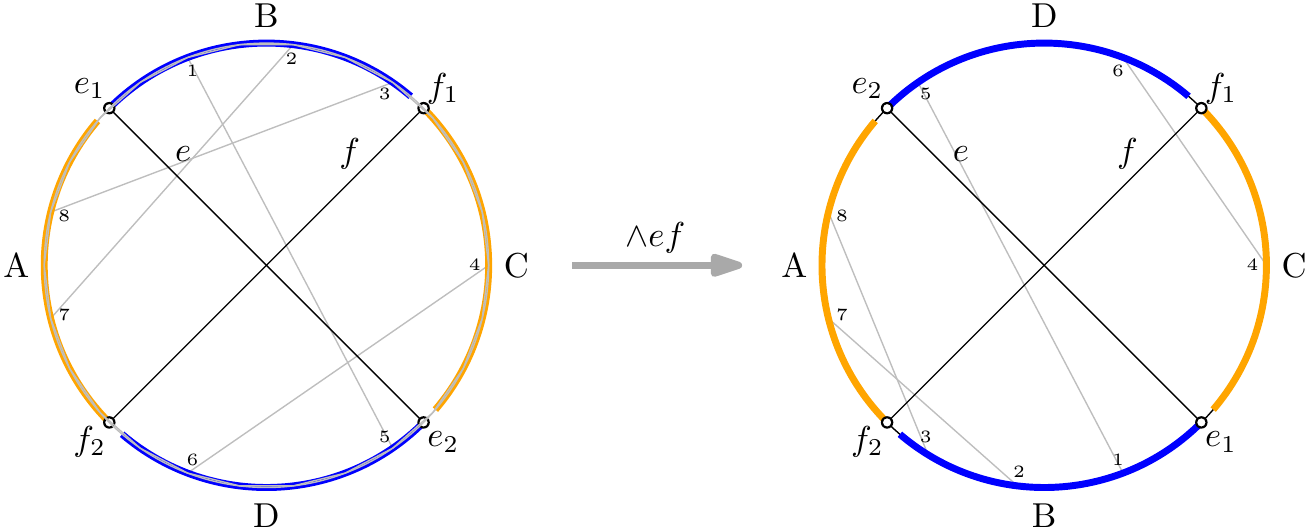}
	\caption{Pivoting $ef$ consists in rotating the parts between the endpoints of the chords $e$ and $f$ as in the picture.}
	\label{fig:pivoting}
\end{figure}

While pivoting is defined on circle graphs, its effect on a chord diagram representation is easy to describe (see \Cref{fig:pivoting}).

On the other hand, contracting an edge in the tree $T$ (or deleting an edge in its complement) results  in the deletion of the corresponding chord from the chord diagram. 
This connection extends to general (loopless simple) graphs, and
this nice interplay partly explains the development of a study of vertex minors and pivot minors, which was built as an analog to Robertson and Seymour's graph minor project.

However, considering graphs instead of maps, and circle graphs instead of maps, we might lose some geometrical and topological aspects. Therefore, we shall keep chord diagrams and maps as primitive objects.

\section{Tours, Quasi-trees, and Map Minors}
\begin{definition}[tour]
Let $\mathcal M=(B,\sigma,\alpha)$ be a general map and let $F$ be a subset of edges of $\mathcal M$. The \emph{tour} of $F$ in $\mathcal M$ is the permutation $\tau$ defined by
\begin{equation}
	\tau(b)=\sigma\alpha_F(b)=\begin{cases}
		\sigma\alpha(b)&\text{if }\underline{b}\in F\\
		\sigma(b)&\text{otherwise}
	\end{cases}
\end{equation}
where, by a slight abuse of notation, we denote by $\alpha_F$ the cutting out of $\alpha$ on the union of all the edges in $F$, that is on the set $\{b\in B\colon \underline b\in F\}$.
\end{definition}
Note that the tour of $F$ in $\mathcal M$ is the same as the tour of the complement of $F$ in the dual general map $\mathcal M^*$.

\begin{figure}[H]
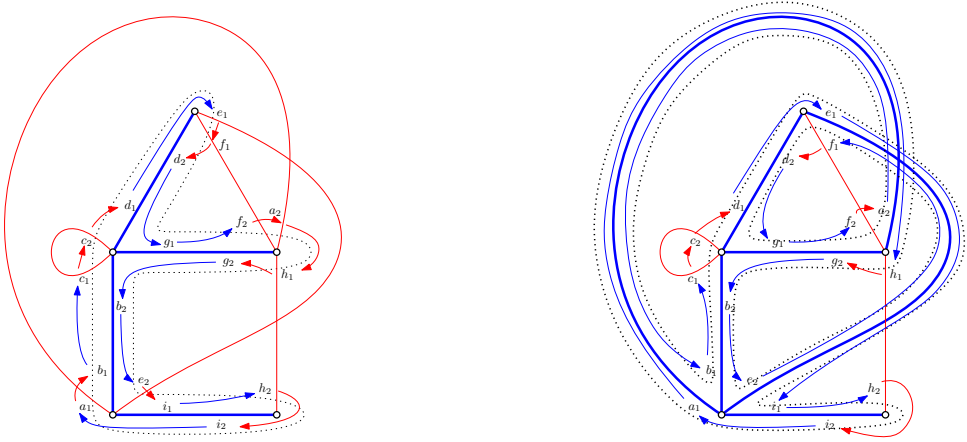

	\centering
	\includegraphics[width=0.4\textwidth]{tour}\hfill
	\includegraphics[width=0.4\textwidth]{tour2}
	\caption{Tour of a spanning tree (on the left) and of a quasi-tree (on the right). (Quasi-)tree edges  are depicted as thick blue  edges and the other edges are depicted as thin red edges.
		The blue arrows correspond to $\sigma\alpha$ and the red ones to $\sigma$.}
	\label{fig:tau}
\end{figure}

When $\mathcal M$ has a single component (i.e. $\mathcal M$ is a map), 
and $F$ is (the edge set of) a spanning tree of $G(\mathcal M)$,  it is well known that  the tour of $F$ in $\mathcal M$ is a cycle (see \Cref{fig:tau}).

We use this property to generalize the notion of spanning tree.

\begin{definition}[quasi-tree]
A  (spanning) \emph{quasi-tree} of a general map $\mathcal M$ is a subset $F$ of edges of $\mathcal M$, whose tour  in $\mathcal M$ is a cycle.
\end{definition}

Note that the standard definition of a quasi-tree is slightly different, as it is (topologically) defined as a map (or, more generally, a ribbon graph) with one boundary component \cite{champanerkar2011quasi}. However, it is not difficult to check that our definition coincides with the notion of a (spanning) quasi-tree in its original sense. The \emph{genus}
of a quasi-tree $S$ of a map $\mathcal M$ is the value $(|S|-(n-1))/2$, where $n$ is the number of vertices of $\mathcal M$.

The existence of a quasi-tree characterizes maps, the same way that the existence of a spanning tree characterizes connected graphs.

\begin{lemma}
	A general map  $\mathcal M=(B,\sigma,\alpha)$  has a single component (i.e. is  a map) if and only if $\mathcal M$ has a quasi-tree.
\end{lemma}	
\begin{proof}
	If $\mathcal M$ is a map, then any spanning tree of $G(\mathcal M)$ is a quasi-tree of $\mathcal M$.
	Conversely, assume that $\mathcal M$ has a quasi-tree $S$. As $\langle\sigma,\alpha\rangle$  acts transitively on each cycle of the tour of $S$ in $\mathcal M$, it acts transitively on $B$. Hence, $\mathcal M$ is a map.
\end{proof}

Quasi-trees are  compatible with map duality, in a way which is similar to trees with matroid duality.
\begin{lemma}
	A subset $S$ of edges is a quasi-tree of a map $\mathcal M$  if and only if its complement is a quasi-tree of the dual map $\mathcal M^\ast$.
\end{lemma}
\begin{proof}
This directly follows from the fact that the tour of $S$ in $\mathcal M$ is the same  as the tour of the complement of $S$ in $\mathcal M^\ast$.
\end{proof}

\begin{figure}[ht]
	\centering
	\includegraphics[width=.75\linewidth]{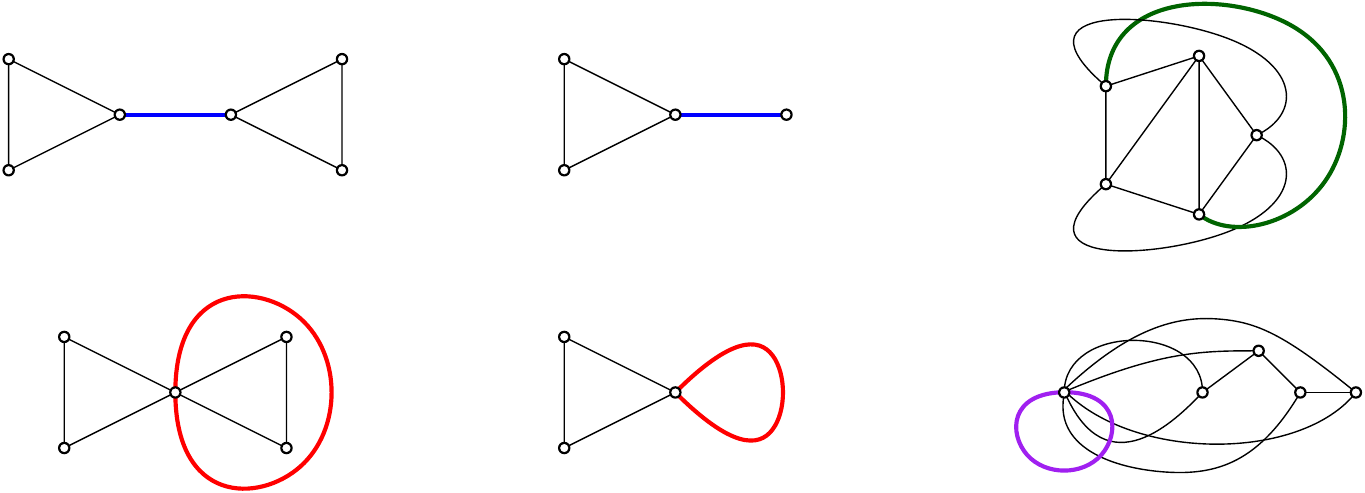}
	\caption{Bridges (in fat blue) and separating loops (in fat red). Bridges belong to all quasi-trees, while separating loops  belong to no quasi-tree. On the right, two dual maps (you can check it\dots), where the purple edge is dual to the green edge.
	On the right top,	the green edge is an example of a non-bridge with two incidences with a same face; on the right bottom, the purple edge is an example of a  non-separating loop.}
	\label{fig:bridge}
\end{figure}

A \emph{bridge} (or \emph{isthmus}) of a map $\mathcal M$ is an edge $e$, such that $G(\mathcal M)-e$ is disconnected, 
while a \emph{separating loop} of $\mathcal M$ is an edge $e$, such that $G(\mathcal M^\ast)-e$ is disconnected. Hence, a separating loop of a map is a bridge of the dual map. 
Remark that bridges (resp. separating loops) are exactly the edges of $\mathcal M$ that belong to all quasi-trees (resp. to no quasi-tree) of $\mathcal M$.
Note that, from the map point of view, a bridge does not always disconnect the map as it can be an edge incident to a degree $1$ vertex. Similarly, a separating loop can be a loop that enclosing a single face (see \Cref{fig:bridge}). 

\begin{definition}[edge deletion]
Let $\mathcal M$ be a map, and let $e$ be an edge of $\mathcal M$ that is not a bridge of $\mathcal M$. The map
 obtained by deleting $e$ in $\mathcal M$ is the map $\mathcal M\del e=(B\setminus e,\sigma_{|B\setminus e},\alpha_{|B\setminus e})$.
%
\end{definition}

\begin{figure}[ht]
	\centering
	\includegraphics[width=\textwidth]{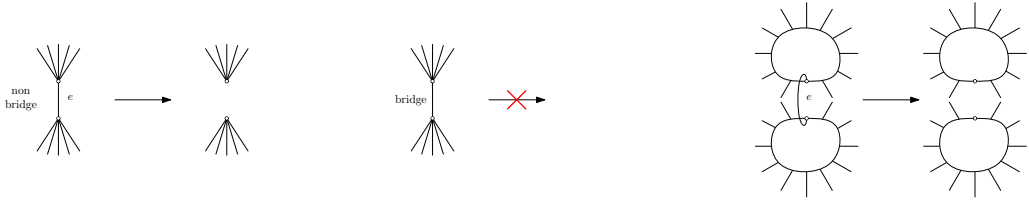}
	\caption{Deleting an edge $e$ of a map $\mathcal M$. The number of edges decreases by $1$, while the number of vertices remains unchanged. Note that the number of faces can increase.}
	\label{fig:deletion}
\end{figure}

We now prove that this definition is valid, meaning that $\mathcal M\del e$ is indeed a map. This will follow from the next lemma applied to any spanning tree $F$ of $G(\mathcal M)$ that avoids $e$,  the transitive action of $\langle\sigma',\alpha'\rangle$  on $B\setminus e$  being witnessed by the existence of a quasi-tree of $\mathcal M\del e$.
Remark that $G(\mathcal M\del e)=G(\mathcal M)\del e$.

\begin{lemma}
	Let $\mathcal M=(B,\sigma,\alpha)$ be a general map, 
let $S\subset B$,  let $e\in B\setminus S$, let $\tau$ be the tour of $S$ in $\mathcal M$, and let $\tau'$ be the tour of $S$ in $\mathcal M\setminus e$. 
Then, $\tau'=\tau_{|B\setminus e}$. 

In particular, $S$ is a quasi-tree of $\mathcal M$ if and only if $S$ is a quasi-tree of $\mathcal M\setminus e$.
\end{lemma}
\begin{proof}
	By definition, $\tau=\sigma\alpha_S$.
 As $e\notin S$ we have $(\alpha_S)_{B\setminus e}=\alpha_S$.
	According to \Cref{lem:prodrest}, we have
	$\tau_{|B\setminus e}=(\sigma\alpha_S)_{|B\setminus e}=\sigma_{|B\setminus e}(\alpha_S)_{|B\setminus e}=\sigma_{|B\setminus e}(\alpha_{|B\setminus e})_S=\tau'$.
\end{proof}

\begin{lemma}
	\label{lem:coloop}
	Let $\mathcal M$ be a map and let $e$ be a bridge of $\mathcal M$. Then $e$ belongs to all the quasi-trees of $\mathcal M$.
\end{lemma}
\begin{proof}
	Assume for contradiction that $e$ is a bridge of $\mathcal M$ and $\mathcal M$ has a quasi-tree $S$ that does not contain $e$.  Then $S\setminus e$ is a quasi-tree of $\mathcal M\del e$, thus $\mathcal M\del e$ is a map. However, it is easily checked that $G(\mathcal M\del e)=G(\mathcal M)\del e$, thus is not connected.
\end{proof}

Note that if we delete from $\mathcal M$ all the edges that are not in a quasi-tree $S$ of $\mathcal M$, we get a map with a single face (which is the tour of $S$ on this map).
Also note that deleting an edge keeps the number of vertices constant.

\begin{definition}[edge contraction]
	Let $\mathcal M$ be a map, and let $e$ be an edge of $\mathcal M$ that is not a separating loop of $\mathcal M$. The map
	obtained by contracting $e$ in $\mathcal M$ is the map $\mathcal M\cont e=(B\setminus e,(\sigma\alpha)_{|B\setminus e}\alpha_{|B\setminus e},\alpha_{|B\setminus e})$.
\end{definition}

\begin{figure}[ht]
	\centering
	\includegraphics[width=\textwidth]{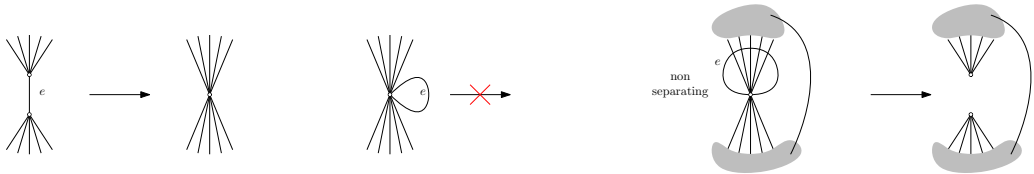}
	\caption{Contracting an edge $e$ of a map $\mathcal M$. Note the special case of a non-separating loop. The number of edges decreases by $1$, while the number of faces remains unchanged. 		}
	\label{fig:contraction}
\end{figure}

It is immediate from the definition that we have 
	\[
	\mathcal M\cont e=(\mathcal M^\ast\del e)^\ast.
	\]	
	As a consequence, if $S$ is a quasi-tree of $\mathcal M$ and $e\in S$ then $S\setminus e$ is a quasi-tree of $\mathcal M\cont e$.

Using the properties we know about quasi-trees and minors, we can easily prove that deletion and contraction commute:
\begin{lemma}
	Let $e\neq f$ be edges of $\mathcal M$. If $\mathcal M\del e\cont f$ is defined then so is $\mathcal M\cont f\del e$ and we have
	$\mathcal M\del e\cont f=\mathcal M\cont f\del e$.
\end{lemma}
\begin{proof}
	Let $T$ be spanning tree of $\mathcal M\del e\cont f$. Then $T$ is a quasi-tree of $\mathcal M\del e$ and $T\cup\{e\}$ is a quasi-tree of $\mathcal M$.
	As $f\notin T\cup\{e\}$ the edge $f$ is contractible in $\mathcal M$ and $T\cup\{e\}$ is a quasi-tree of $\mathcal M\cont f$. As $e\in T\cup\{e\}$, the edge $e$ is erasable in $\mathcal M\cont f$ and $T$ is a quasi-tree of $\mathcal M\cont f\del e$. As the tours of $T$ in $\mathcal M\del e\cont f$ and $\mathcal M\cont f\del e$ are the same, the two maps are the same.
\end{proof}

Note that a similar proof can be used to prove that deletions commute and that contractions commute.
From the properties of quasi-trees with respect to deletion and contraction,  we deduce the following.

\begin{lemma}
Let  $\varsigma(\mathcal M)$ be the number of quasi-trees of $\mathcal M$.
Then $\varsigma$ satisfies the following contraction/deletion formula
	\[
	\varsigma(\mathcal M)=
	\begin{cases}
		\varsigma(\mathcal M\cont e)&\text{if $e$ is a bridge}\\
		\varsigma(\mathcal M\del e)&\text{if $e$ is a separating loop}\\
		\varsigma(\mathcal M\cont e)+\varsigma(\mathcal M\del e)&\text{otherwise}
	\end{cases}
	\]
\end{lemma}	

\begin{figure}[ht]
	\centering
	\includegraphics[width=.7\textwidth]{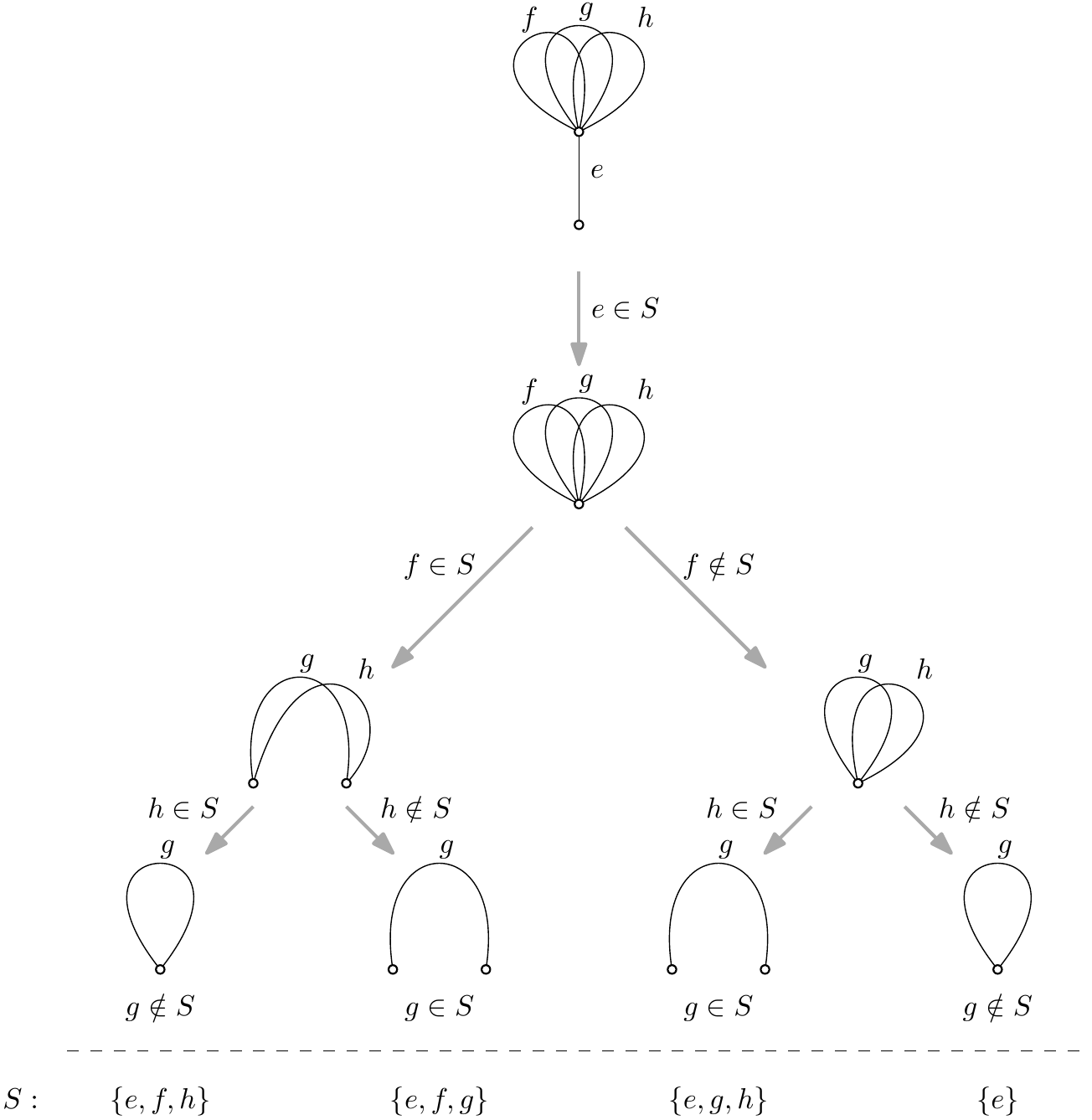}
	\caption{A map with $4$ quasi-trees}
	\label{fig:delcont}
\end{figure}

Notice that a map $\mathcal M$ defines two graphs, namely $G(\mathcal M)$ and $G(\mathcal M^\ast)$. We noticed that edge deletion corresponds to a deletion in $G(\mathcal M)$, thus edge contraction corresponds to a deletion in $G(\mathcal M^\ast)$. 

For more on this subject and on the extension to maps (and, more generally, to ribbon graphs) of the Tutte polynomial, we refer the reader to \cite{champanerkar2011quasi} and to the survey~\cite{chmutov2017topological}.

\section{Bicolored chord diagrams and the $\Delta$-matroid of quasi-trees}

Given a cyclic permutation $\tau$, we define the \emph{interlace chord diagram} $\Lambda(\tau)$ of $\tau$ as the chord diagram obtained by putting on a circle the flags in $B$ in $\tau$-order, the chords linking the pairs of flags in a same edge.
Given a quasi-tree $S$ of a map $\mathcal M$ with tour $\tau$, we define $\widetilde{\Lambda}(\mathcal M, S)$ as the chord diagram $\Lambda(\tau)$ with chords in $S$ colored $1$ and those not in $S$ colored $2$ (see \Cref{fig:tauCD}). Every chord diagram is a representation of a \emph{circle graph}, whose vertices are the chords and whose edges are the pairs of intersecting chords. It will be convient to denote $\widetilde{\Lambda}(\mathcal M, S)$  both the chord diagram and its associated circle graph.

\begin{figure}[ht]
	\centering
	\includegraphics[width=0.4\textwidth]{tour2}\hfill
	\includegraphics[width=0.4\textwidth]{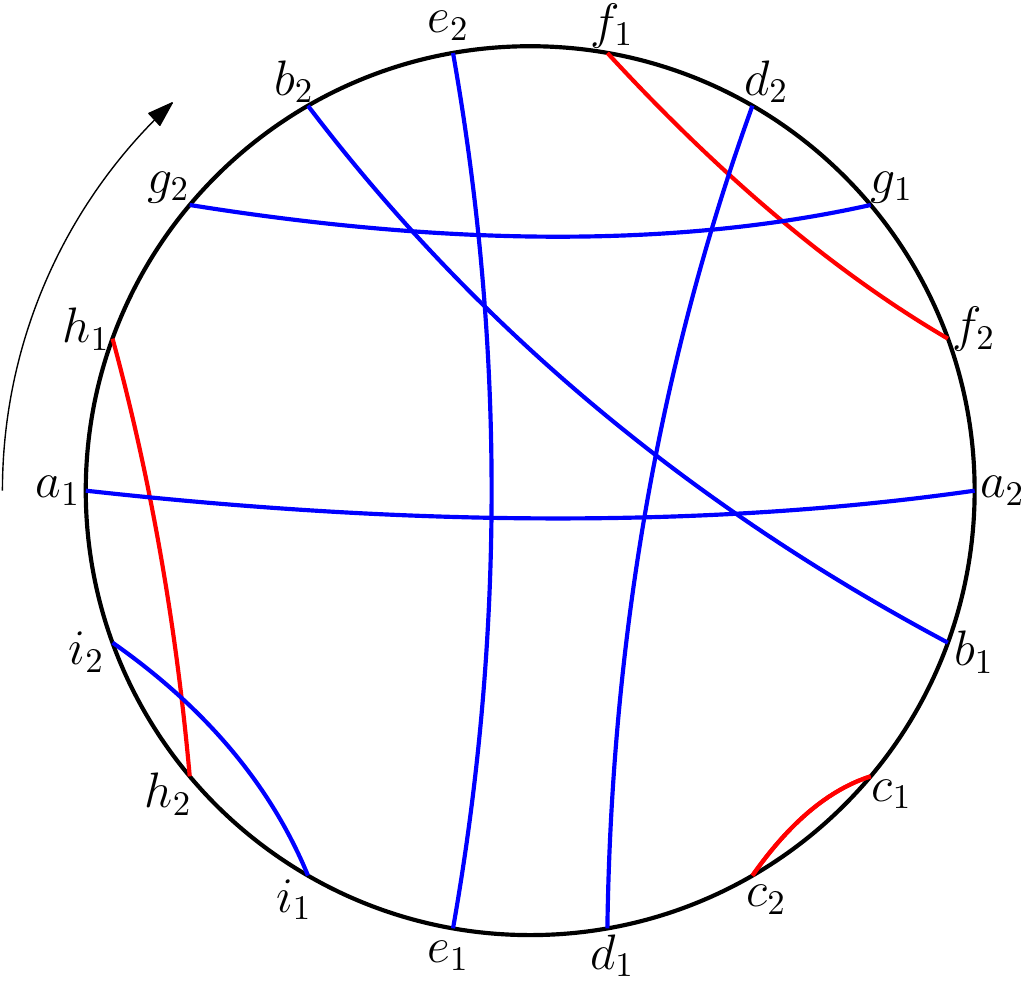}
	\caption{Tour of a quasi-tree and the corresponding bicolored chord diagram (color $1$ is blue, color $2$ is red).}
	\label{fig:tauCD}
\end{figure}

We further define as $I(\mathcal M,S)$ the bipartite subgraph of $\widetilde{\Lambda}(\mathcal M,S)$ induced by $S$ and its complement. In other words, two edges are adjacent in 
$I(\mathcal M,S)$  if they are adjacent in $\widetilde{\Lambda}(\mathcal M,S)$,  one is in $S$ and not the other. The property mentioned in the planar case extends in the general case:
\begin{fact}
	\label{fact:YZ}
	Let $\mathcal M$ be a map and let $S$ be a spanning tree of $\mathcal M$. Then 
\begin{itemize}
	\item for each $e\in S$, the neighbors of $e$ in $I(\mathcal M,S)$ are the elements of $\omega(e)\setminus\{e\}$;
	\item for each $f\notin S$, the neighbors of $f$ in $I(\mathcal M,S)$ are the elements of $\gamma(f)\setminus\{f\}$.
\end{itemize}
\end{fact}

Note that a direct consequence of \Cref{fact:YZ} is that $I(\mathcal M,S)$ (hence $\widetilde{\Lambda}(\mathcal M,S)$) is connected if $\mathcal M$ is $2$-connected.

It is well known that two spanning trees $S$ and $S'$ differ by two elements, i.e. $S'=S\bigtriangleup\{e,f\}$  if and only if $e$ and $f$ are adjacent in $I(\mathcal M,S)$, and that  we have $I(\mathcal M,S')=I(\mathcal M,S)\wedge ef$. This property extends to the chord diagrams of quasi-trees (considering $\widetilde{\Lambda}(\mathcal M,S)$ instead of $I(\mathcal M,S)$).

\begin{lemma}
	\label{lem:pivoting_qt}
	Let $S$ be a quasi-tree of a map $\mathcal M$ and let $e,f$ be distinct edges of $\mathcal M$ that are adjacent in $\widetilde{\Lambda}(\mathcal M, S)$. Then $S'=S\bigtriangleup\{e,f\}$ is a quasi-tree of $\mathcal M$ and $\widetilde{\Lambda}(\mathcal M, S')$ is obtained from $\widetilde{\Lambda}(\mathcal M, S)\wedge ef$ by flipping the colors of $e$ and $f$.
	
	Conversely, if $S'=S\bigtriangleup\{e,f\}$ is a quasi-tree of $\mathcal M$, then $e$ and $f$ are adjacent in $\widetilde{\Lambda}(\mathcal M, S)$. 
\end{lemma}
\begin{proof}
	Let $\tau$ be the tour of $S$ in $\mathcal M$, and let $e_1,f_1,e_2,f_2$ be the flags in $e$ and $f$ (following the $\tau$-order). The cycle $\tau$ then rewrites as $(w_1,e_1,w_2,f_1,w_3, e_2, w_4, f_2)$, where $w_1,w_2,w_3$, and $w_4$ are sequences of flags.
	It is easily checked that the tour $\tau'$ of $S'$ is $(w_1, e_2, w_4, f_1,w_3,e_1,w_2)$. In particular, $\tau'$ is a cycle thus $S'$ is a quasi-tree of $\mathcal M$.
	
	Assume that $e$ and $f$ are not djacent in $\widetilde{\Lambda}(\mathcal M, S)$. Let $e_1,e_2,f_1,f_2$ be the flags in $e$ and $f$ (following the $\tau$-order). The cycle $\tau$ then rewrites as $(w_1,e_1,w_2,e_2,w_3, f_1, w_4, f_2)$, where $w_1,w_2,w_3$, and $w_4$ are sequences of flags.
Then, it is easily checked that the tour $\tau'$ of $S'$ has $3$ cycles, namely $(w_1, e_1, w_3, f_1)$, $(w_2,e_2)$ and $(w_4,f_2)$. Thus, $S'$ is not a quasi-tree of $\mathcal M$.
\end{proof}

\begin{figure}[H]
	\centering
	\includegraphics[width=0.4\textwidth]{tour2CD}\hfill
	\includegraphics[width=0.4\textwidth]{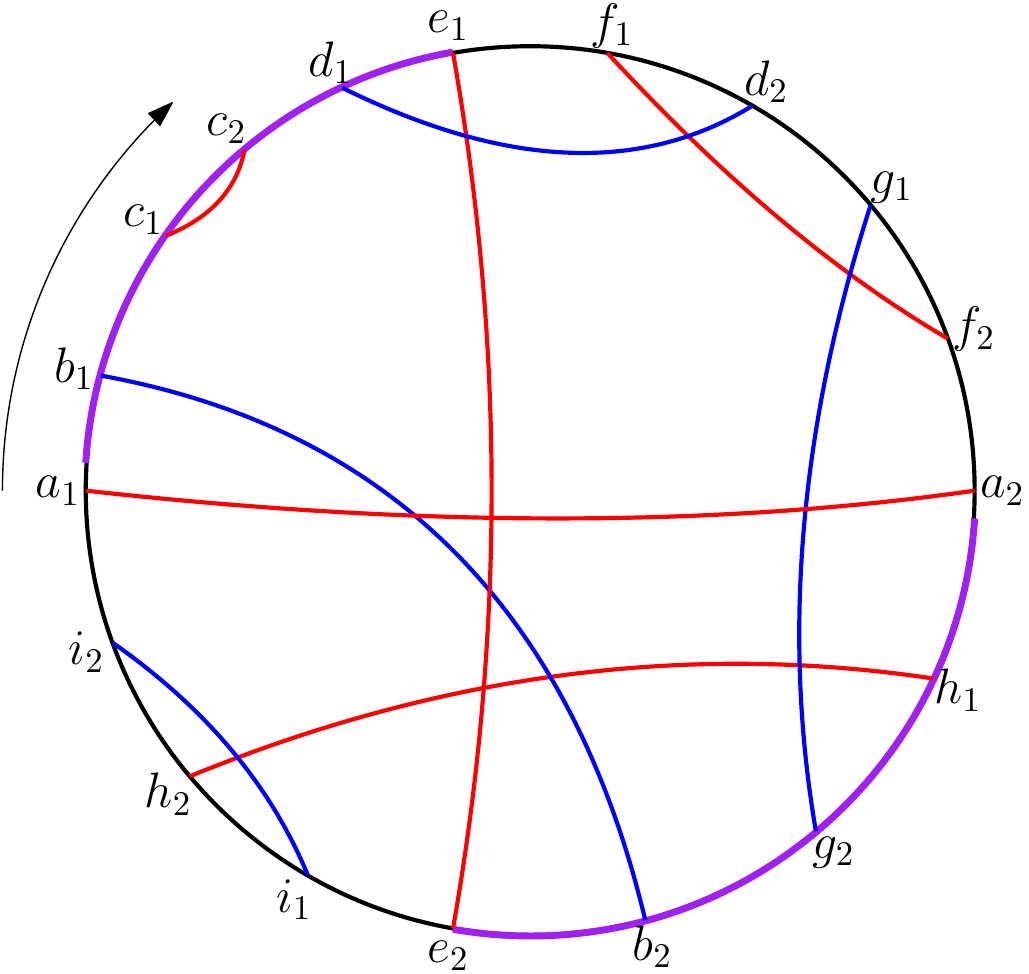}
	\caption{Pivoting $ae$ transforms the chord diagram of the quasi-tree $\{a,b,d,e,g,i\}$ into the one of the spanning-tree $\{b,d,g,i\}$.}
	\label{fig:tauCD}
\end{figure}

\begin{lemma}
	\label{lem:YZ}
	For every quasi-tree $S$ of $\mathcal M$, there exist quasi-trees $S_0$ and $S_1$ of $\mathcal M$ with
	$S_0\subseteq S\subseteq S_1$, $S_0$ is a spanning tree of $\mathcal M$, and $S_1$ is the complement of a spanning tree of $\mathcal M^\ast$.
\end{lemma}
\begin{proof}
	By duality, it is sufficient to prove the existence of $S_0$. Assume that $S$ contains two edges $e,f$ that are adjacent in $\widetilde{\Lambda}(S)$. Then $S'=S\bigtriangleup\{e,f\}\subset S$ is a quasi-tree of $\mathcal M$. So, assume that $S$ does not contain two edges $e,f$ that are adjacent in $\widetilde{\Lambda}(\mathcal M, S)$. By deleting all the edges not in $S$, we get a map $\mathcal M'$, where $G(\mathcal M')$ is a spanning subgraph of $G(\mathcal M)$ induced by the edges in $S$. As this map has a single face and no interlaced edges, we deduce that $\mathcal M'$ is a tree, hence $S$ is a spanning tree of $\mathcal M$.
\end{proof}

The next lemma is an easy consequence of our study of minors.
\begin{lemma}
	Let $S$ be a quasi-tree of a map $\mathcal M$. 
	\begin{itemize}
\item 	If $e\in S$, then $\widetilde{\Lambda}(\mathcal M\cont e, S\setminus e)=\widetilde{\Lambda}(\mathcal M, S)- e$;
\item If $e\notin S$, then $\widetilde{\Lambda}(\mathcal M\del e, S)=\widetilde{\Lambda}(\mathcal M, S)- e$.
	\end{itemize}
\end{lemma}


We now state an analog of the edge exchange property of spanning trees, which implies that the set of all quasi-trees of a map forms an  \emph{even $\Delta$-matroid} \cite{Bouchet1987}. 
\begin{lemma}
	\label{lem:exch}
	Let $S_1\neq S_2$ be two quasi-trees of a map $\mathcal M$. Then, for every $e\in S_1\bigtriangleup S_2$ there exists $f\in S_1\bigtriangleup S_2$ with $f\neq e$, such that $S_1\bigtriangleup\{e,f\}$ is a quasi-tree of $\mathcal M$.
\end{lemma}
\begin{proof}
	By contracting all the edges in $S_1\cap S_2$ and deleting those in the complement of $S_1\cup S_2$, we reduce to the case where $(S_1,S_2)$ is a partition of the edge set. Note that  $\mathcal M$ has no bridges and no separating loops.
	Let $\tau_1$ be the tour of $S_1$ in $\mathcal M$. If there exists $f\neq e$ such that $ef$ is an edge of $\widetilde{\Lambda}(\mathcal{M},S_1)$ we are done. We prove by contradiction that no other case can occur.

		Assume $e\in S_1$.   Let $\tau_1=(w_1,e_1,w_2,e_2)$. 
		Note that $w_1$ and $w_2$ are not empty. Let $\mathcal M\del e=(B\setminus e,\sigma',\alpha')$.
		Let $b\in w_2$. By assumption, $\alpha(b)\in w_2$. 
		Assume $b\notin S_1$. If $b$ is not the last flag of $w_2$ then $\sigma'(b)=\sigma(b)\in w_2$. Otherwise, $\sigma'(b)=\sigma^2(b)=\sigma(e_2)=\sigma\alpha(e_1)\in w_2$.
		Assume $b\in S_1$ (hence $\alpha(b)\in S_1$). If $\alpha(b)$ is not the last flag of $w_2$ (i.e. $\sigma(b)\neq e_2$) then $\sigma'(b)=\sigma\alpha(\alpha(b))\in w_2$. Otherwise, $\sigma'(b)=\sigma^2(b)=\sigma(e_2)=\sigma\alpha(e_1)\in w_2$. In all cases, $\sigma'(b)\in w_2$.
		Thus, we get that $w_2$ is closed under the action of $\langle \sigma',\alpha'\rangle$, contradicting the hypothesis that $e$ is not a bridge (i.e. that $\mathcal M\del e$ is a map).
		
		Assume $e\notin S_1$.   Let $\tau_1=(w_1,e_1,w_2,e_2)$. 
		Note that $w_1$ and $w_2$ are not empty. Let $\mathcal M\cont e=(B\setminus e,\sigma',\alpha')$.
		Let $b\in w_2$. By assumption, $\alpha(b)\in w_2$. 
		Assume $b\in S_1$. If $b$ is not the last flag of $w_2$ then $\sigma'\alpha'(b)=\sigma\alpha(b)\in w_2$. Otherwise, $\sigma'\alpha'(b)=\sigma\alpha\sigma\alpha(b)=\sigma\alpha(e_2)=\sigma(e_1)\in w_2$.
		Assume $b\notin S_1$. If $\alpha(b)$ is not the last flag of $w_2$ (i.e. $\sigma(b)\neq e_2$) then $\sigma'\alpha'(b)=\sigma\alpha(b)\in w_2$. 
		Otherwise, $\sigma'\alpha'(b)=\sigma\alpha\sigma\alpha(b)=\sigma\alpha(e_2)=\sigma(e_1)\in w_2$. In all cases, $\sigma'\alpha'(b)\in w_2$.
		Thus, we get that $w_2$ is closed under the action of $\langle \sigma'\alpha',\alpha'\rangle$, contradicting the hypothesis that $e$ is not a separating loop (i.e. that $\mathcal M\cont e$ is a map).
\end{proof}

\begin{corollary}
	\label{cor:qt_connect}
	Let $S_1\neq S_2$ be two quasi-trees of a map $\mathcal M$.  Then there exists a sequence $(\{e_1,f_1\},\dots,\{e_k,f_k\})$ of pairs of edges such that 
	$S_2=S_1\bigtriangleup\{e_1,f_1\}\bigtriangleup\dots\bigtriangleup\{e_k,f_k\}$ and, 
	for all $1\leq i<k$, $S_1\bigtriangleup\{e_1,f_1\}\bigtriangleup\dots\bigtriangleup\{e_i,f_i\}$ is a quasi-tree of $\mathcal M$.
\end{corollary}

\begin{remark}
	\label{rem:loop_coloop}
	In this setting,  the bridges (resp. the separating loops) of $\mathcal M$ are the \emph{coloops} (resp. the \emph{loops}) of the $\Delta$-matroid of the quasi-trees of $\mathcal M$, meaning that bridges are exactly those edges that belong to all quasi-trees (\Cref{lem:coloop}) and, dually, separating loops are exactly those edges that belong to no quasi-tree.
	
	Also, as a consequence of \Cref{cor:qt_connect,lem:pivoting_qt}, the isolated vertices of $\widetilde{\Lambda}(\mathcal{M},S)$ are exactly the bridges and the separating loops of $\mathcal M$.
\end{remark}

\section{Representation of maps using quasi-trees}
We extend the bijection between bipartite circle graphs and planar graphs shown \Cref{fig:bipcirc} to a bijection between bicolored chord diagrams
and pairs $(\mathcal M,S)$, where $\mathcal M$ is a map and $S$ is a quasi-tree of $\mathcal M$. 

We now present how to derive a representation of a map  $\mathcal M$ from the choice of a quasi-tree $S$ of $\mathcal M$. Let $g$ be the genus of $\mathcal M$ and let $g_S$ be the genus of $S$.

\begin{figure}[ht]
	\centering
	\includegraphics[width=.75\textwidth]{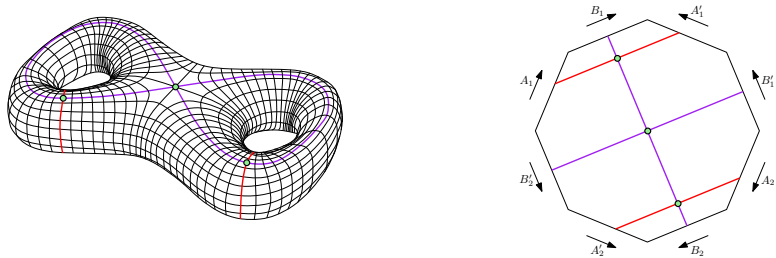}
\caption{Representation of a double torus by an octogon.}
\label{fig:poly}
\end{figure}

Recall that an orientable surface of genus $g>1$ can be represented using 
a polygon with $4g$ sides (the \emph{Fricke canonical polygon} or  \emph{fundamental polygon} of the surface), where the sides are matched in distinct pairs with opposite orientation (see \Cref{fig:poly} and, e.g. \cite{moharthom}).

\begin{figure}[ht]
	\centering
	\begin{minipage}{.3\textwidth}
		\centering \includegraphics[width=.8\textwidth]{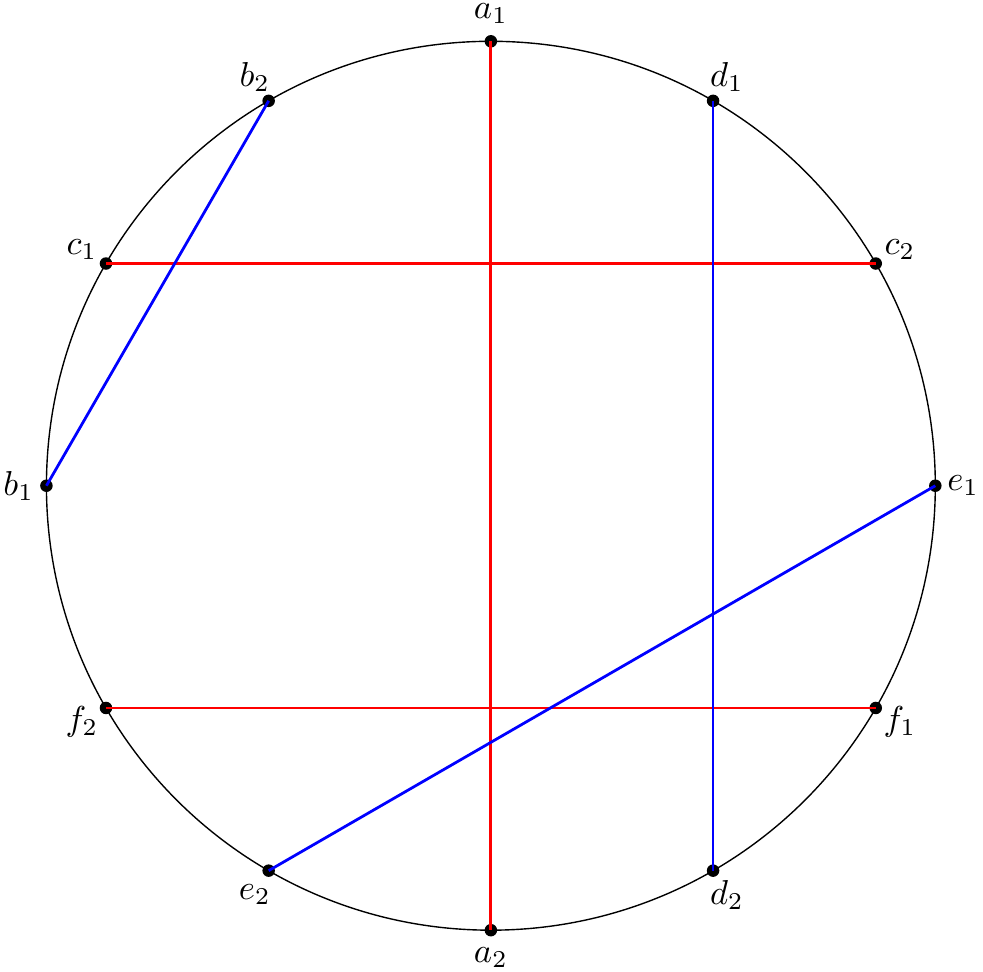}
	\end{minipage}
	\hfill
	\begin{minipage}{.3\textwidth}
		\includegraphics[width=\textwidth]{BCD4-3}
	\end{minipage}
	\hfill
	\begin{minipage}{.3\textwidth}
		\includegraphics[width=\textwidth]{BCD4-2}
	\end{minipage}
	
	\begin{minipage}{.3\textwidth}
		\includegraphics[width=\textwidth]{BCD4-1}
	\end{minipage}
	\hfill
	\begin{minipage}{.3\textwidth}
		\includegraphics[width=\textwidth]{BCD4}
	\end{minipage}
	\hfill
	\begin{minipage}{.3\textwidth}
		\includegraphics[width=\textwidth]{BCD4-0}
	\end{minipage}
	
	\caption{A bicolored chord diagram  (top left), a representation of the double torus using two squares (center), and a drawing of the chord diagram without crossings (top right).
		Dualization of the interior map (bottom left), connection to the endpoints of the chords not in $S$ (bottom center), and final drawing of the map (and quasi-tree) associated to the bicolored chord diagram (bottom right).}
	\label{fig:drawing}
\end{figure}

In our representation of maps, we  use up to two polygons.
Precisely, we fix a circle $\Gamma$, a polygonal $\Pi_i$ inside $\Gamma$ with $4g_S$ sides and a polygon $\Pi_e$ enclosing $\Gamma$ with $4(g-g_S)$ sides. (If $g_s=0$ we don't need $\Pi_i$ and if $g=g_s$ we don't need $\Pi_e$.) Each of these two polygons has his sides matched by distinct pairs, as in Fricke representation.
We put on $\Gamma$ the points of the chord diagram $\widetilde{\Lambda}(\mathcal M,S)$ and  draw (without crossings, using the fundamental polygons) the chords of $S$ inside $\Gamma$ and the other chords outside $\Gamma$. These chords are drawn in such a way that they do not cross $\Gamma$ (see \Cref{fig:drawing}).


If one contracts $\Gamma$ into a single point, our drawing defines two maps: the exterior map with edges ``outside''  $\Gamma$, which has genus $g_S$, and the interior map with edges ``inside'' $\Gamma$, which has genus $g-g_S$. Note that each of these maps has a single vertex (corresponding to $\Gamma$). We now draw the dual of the interior map in our representation, inside $\Gamma$. Note that the obtained map has a number of vertices equal to $|S|-2g_S+1$.
Then, we connect each point of $\Gamma$ corresponding an endpoint of a chord not in $S$ to the vertex that can be reached without crossing any other edge.

\section{Rooted maps, Ordered matchings, and Double occurence words}

Recall that rooting a map $\mathcal M=(B,\sigma,\alpha)$ consists in selecting a root flag $b_\bullet\in B$. Note that rooting a map kills all the automorphisms of the map, as every flag can be written as $\mu(b_\bullet)$, for some $\mu\in\langle\sigma,\alpha\rangle$.  It will be convenient to denote $\mathcal M_\bullet$ the map $\mathcal M$ rooted at $b_\bullet$.

\begin{figure}[ht]
	\centering
	\includegraphics[width=.4\textwidth]{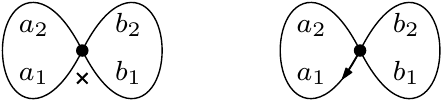}
	\caption{A rooted map $\mathcal M_\bullet$. The root flag is either indicated using a cross placed just (anti-clockwise) before the root flag (here $a_1$) at the incident vertex (on the left), or by orienting the corresponding edge from the vertex incident to the root flag (on the right).}
	\label{fig:rooted_map}
\end{figure}

Similarly, rooting a chord diagram consists in selecting a root chord endpoint.
Rooting $\widetilde{\Lambda}(\mathcal M,S)$ and reading all the flags from the root, while traversing the circle clockwise, we obtain a linear order on $B$ or, equivalently, a \emph{single occurrence word}, that is a word whose letters are the flags, which contains every flag exactly once. 
This way, we can identify rooted bicolored chord diagrams and bicolored ordered matchings. For example, the bicolored chord diagram depicted on the left of \Cref{fig:tauCD} and rooted at $a_1$ corresponds to the bicolored ordered matching

\[
\xy 
\POS(0,0)*{}*\cir<2.5pt>{}="a1", 
(5,0)*{}*\cir<2.5pt>{}="h1",
(10,0)*{}*\cir<2.5pt>{}="g2",
(15,0)*{}*\cir<2.5pt>{}="b2",
(20,0)*{}*\cir<2.5pt>{}="e2",
(25,0)*{}*\cir<2.5pt>{}="f1",
(30,0)*{}*\cir<2.5pt>{}="d2",
(35,0)*{}*\cir<2.5pt>{}="g1",
(40,0)*{}*\cir<2.5pt>{}="f2",
(45,0)*{}*\cir<2.5pt>{}="a2",
(50,0)*{}*\cir<2.5pt>{}="b1", 
(55,0)*{}*\cir<2.5pt>{}="c1",
(60,0)*{}*\cir<2.5pt>{}="c2",
(65,0)*{}*\cir<2.5pt>{}="d1",
(70,0)*{}*\cir<2.5pt>{}="e1",
(75,0)*{}*\cir<2.5pt>{}="i1",
(80,0)*{}*\cir<2.5pt>{}="h2",
(85,0)*{}*\cir<2.5pt>{}="i2",
\POS"a1"\ar@/^5ex/@{-}@[blue] "a2",
\POS"b1"\ar@/_5ex/@{-}@[red] "b2",
\POS"c1"\ar@/^1ex/@{-}@[red] "c2",
\POS"d1"\ar@/_3ex/@{-}@[blue] "d2"
\POS"e1"\ar@/_5ex/@{-}@[blue] "e2"
\POS"f1"\ar@/^2ex/@{-}@[red] "f2",
\POS"g1"\ar@/_3ex/@{-}@[blue] "g2",
\POS"h1"\ar@/^7ex/@{-}@[red] "h2",
\POS"i1"\ar@/^2ex/@{-}@[blue] "i2"
\endxy
\]

\begin{lemma}
	The following objects are in bijection:
\begin{itemize}[-]
	\item rooted bicolored chord diagrams with $m$ chords, 
	\item bicolored ordered matchings of size $2m$,
	\item pairs $(\mathcal M_\bullet,S)$, where $\mathcal M_\bullet$ is a rooted map with $m$ 
	edges and $S$ is a quasi-tree of $\mathcal M_\bullet$.
\end{itemize}
\end{lemma}	
\begin{proof}
	The chord diagram is defined by $\tau$, which can be computed from $\mathcal M$ and $S$. The bicoloration is given by $S$. 
	Conversely, from the bicolored chord diagram, we can define $\alpha$ (by the chords) and $\sigma$ (from $\tau$, $\alpha$ and the bicoloration). 
\end{proof}

\begin{corollary}
The sum of $\varsigma(\mathcal M_\bullet)$ over all rooted  maps $\mathcal M_\bullet$ with $m$ edges is $\frac{(2m)!}{m!}$.
\end{corollary}
\begin{proof}
This sum is $2^m$ times the number of ordered matchings of size $2m$.
\end{proof}

The  \emph{pivoting class} of a chord diagram $\Lambda$ is the class of all chord diagrams that can be reached from $\Lambda$ by a sequence of pivoting operations.

\begin{lemma}
	The number of distinct pivoting classes of rooted chord diagrams with $m$ chords equals the number of rooted maps with $m$ edges, and	the number $\varsigma(\mathcal M)$ of quasi-trees of a rooted map $\mathcal M_\bullet$ is the size of the pivoting class associated to $\mathcal M_\bullet$.
\end{lemma}
\begin{figure}[ht]
	\centering
	\includegraphics[width=\textwidth]{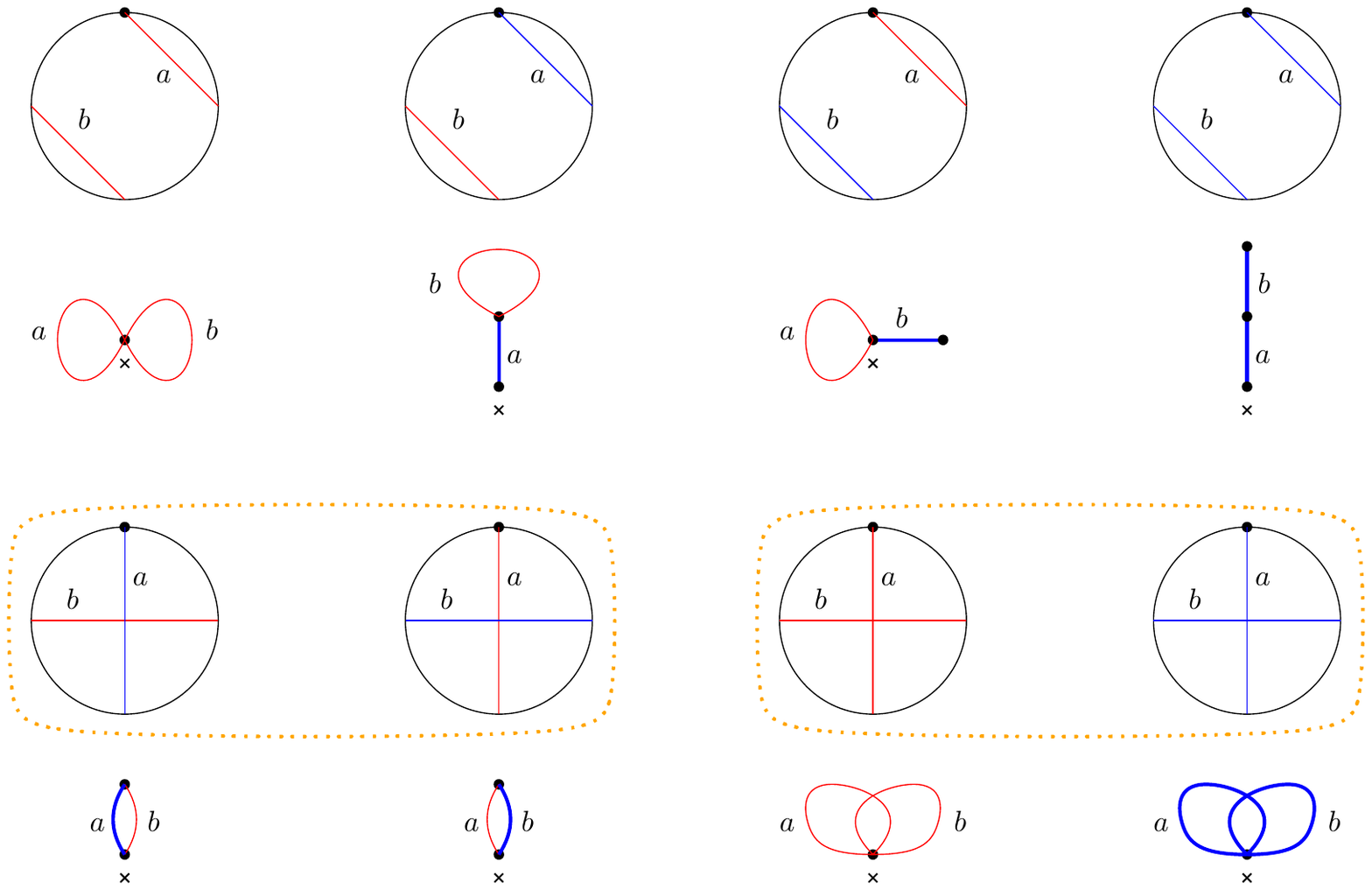}
	\caption{Bijections in the case $m=2$. There are $6$ rooted maps and $8$ pairs $(\mathcal M_\bullet,S)$, where $\mathcal M_\bullet$ is a rooted map and $S$ is a quasi-tree of $\mathcal M$.
	The set $S$ is in blue and doted curves delineate the pivoting classes.}
\end{figure}

 Applying the transformation $b\mapsto \underline{b}$, which maps each flag to the edge it belongs to, we transform the single occurrence word $w$ defined by the tour of a quasi-tree into a word, whose letters are edges, which uses every edge exactly twice. We call such a word a \emph{double occurrence word}. When using the above labeling scheme, a tour of a rooted map is fully determined by the associated double occurrence word. A natural question is whether we can use a standard quasi-tree, which could also be recovered from the mere double occurrence word. In full generality, it is obvious that the answer is negative, as the double occurrence word $e,e$ corresponds to two distinct maps, one with a single edge linking two vertices (with quasi-tree $\{e\}$) and one with a single loop attached to a vertex (with empty quasi-tree). However, if we restrict ourselves to loopless maps, we shall prove that some particular quasi-tree (actually, some specific spanning tree)  can be recovered from its associated double occurrence word and that, moreover, the admissible double occurrence words can be characterized.

\begin{figure}[ht]
	\centering
	\includegraphics[width=\textwidth]{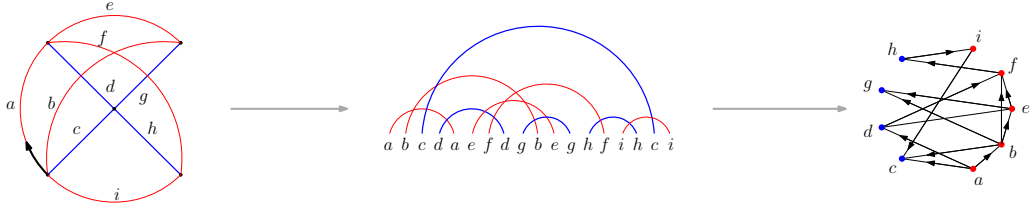}
	\caption{From a spanning tree of a rooted map (on the left) to a directed circle graph (on the right), via a double occurrence word (in the middle).}
	\label{fig:dir_interlace}
\end{figure}

Note that  the circle graph defined by a double occurrence word $w$ gets a natural orientation: the edge $ef$ of the circle graph is oriented from $e$ to $f$ if $e\,f\,e\,f$ is a subword of $w$. We denote by $\vec{\Lambda}(w)$ the oriented circle graph defined by the single occurrence word $w$ (see \Cref{fig:dir_interlace}).

\section{Depth-First Search trees}

The \emph{depth-first search} (DFS) is traversal procedure of a connected graph $G$ that produces a rooted spanning tree $T$ of $G$ with the property that the two incidences of every edge of $G$ lie in a leaf-to-root path of $T$ (\emph{Tr\'emaux property}). It will be convenient to describe this procedure on a rooted combinatorial map $(\mathcal M,b_0)$. The DFS starts at  the vertex incident to $b_0$, which is the \emph{root} of $T$, and explores the graph as far as possible before backtracking. In order to formalize the DFS, we introduce some terminology and notations. 
It is convenient to introduce a dummy flag $b^0$, inserted just before $b_0$ (i.e. $\sigma(b^0)=b_0$).
The current position of the search is a flag $b$, initially set to $b_0$; We denote by $M_V$ (resp $M_B$) the set of \emph{visited vertices} (resp. of \emph{visited flags}). Initially, $M_B=\{b^0\}$ and $M_V=\{r\}$, where $r$ is the vertex incident to $b_0$. To each non-root visited vertex $v\in M_V\setminus\{r\}$ will be associated a \emph{discovery flag} $b_v$, which will be the flag incident to $v$  in the edge linking $v$ to its parent in the constructed spanning tree (or $b^0$ if $v$ is the root). The DFS can be described as the repetition of the following two steps: 
\begin{itemize}
	\item  we look for a non-visited flag incident to the current vertex. If such a flag exists, this becomes the current flag $b$ (which is now visited and consequently added to  $M_B$); otherwise, if the current vertex is  not the root, we backtrack, meaning that the current flag becomes the discovery flag of $v$  (which is now visited and consequently added to  $M_B$); otherwise (i.e. if all the flags incident to the current vertex have been visited and the current vertex is the root), the procedure ends.
	\item if the opposite flag $b'=\alpha(b)$ of the current flag $b$ is not incident to a visited vertex, then we traverse the edge, meaning that the current vertex $v$ is now the vertex incident to $b'$ (which is now visited and consequently added to  $M_V$), the discovery flag of this vertex is set to $b'$, and the current flag is $b'$ (which is now visited and consequently added to  $M_B$).
\end{itemize}
Notice that in this procedure we do not precise how an non-visited flag is chosen when multiple choices are possible. In this sense, the DFS is not (formally speaking) a fully deterministic procedure.

When dealing with combinatorial map, there are two obvious ways of making the above algorithm fully deterministic, by always selecting the first possible choice or by always selecting the last possible choice.
This way, we construct two special DFS-tree, the \emph{Early DFS-tree} and the \emph{Late DFS-tree} of the rooted map.
Note that by reversing the orientation of the map, that is by considering the map $(B,\sigma^{-1},\alpha)$ instead of $(B,\sigma,\alpha)$, the Early DFS-tree (resp. the Late DFS-tree) becomes the Late DFS-tree (resp. the Early DFS-tree).

Given a rooted map $\mathcal M_\bullet$ and a spanning tree $S$ of $\mathcal M_\bullet$, it is easily checked whether $S$ has the Tr\'emaux property, i.e. whether $S$ is a DFS-tree.

\begin{lemma}
	\label{lem:DFS}
	Let $S$ be a spanning tree of a rooted map $\mathcal M_\bullet$, and let $\underline{w}$ be the associated double occurrence word. 
	
	Then $S$ has the Tr\'emaux property if and only if no $f$ has both (in $\vec{\Lambda}(w)$) an in-neighbor and an out-neighbor in $S$. In other words $S$ has the Tr\'emaux property if and only $w$ does not contain the pattern $e_1\,f\,e_1\,e_2\,f\,e_2$, with $e_1,e_2\in S$.
\end{lemma}
\begin{proof}
	First note that if an edge $f$ has a neighbor in $S$, then $f\notin S$.
	Assume $S$ is a DFS-tree.
	Contracting tree edges and deleting cotree edges preserves the property of being a DFS-tree. This allows us to reduce to the case where $M$ is a cycle. Then either $w=f\,e_1\dots e_k\,f\,e_k\dots e_1$ or  $w=e_1\dots e_k\,f\,e_k\dots e_1\,f$. 
	
	Conversely, assume that for every $f\notin S$ either all the $S$-neighbors of $f$ in $\vec\Lambda(w)$ are in-neighbors of $f$ or all the $S$-neighbors of $f$ in $\vec\Lambda(w)$
	are out-neighbors of $f$ and assume for contradiction that $E\setminus S$ contains a transversal edge $f=uv$. Let $x$ be the least common ancestor of $u$ and $v$ in $S$, and let $e_1,\dots,e_k$ (resp. $e_1',\dots,e_\ell'$) be the path linking $x$ and $u$ (resp. $x$ and $v$) in $S$. Again, we can   reduce to case where $M$ is the cycle $\gamma(f)$. Then we get $w=e_1\dots e_k\,f\,e_k\dots e_1\,e_1'\dots,e_\ell'\,f\,e_\ell'\dots,e_1$ (up to exchange between the paths $e_1\dots e_k$ and $e_1'\dots e_\ell'$). Hence,  $f$ has both in-neighbors and out-neighbors in $S$, contradicting our hypothesis.
\end{proof}

\begin{lemma}
	\label{lem:extDFS}
	Let $S$ be a spanning tree of a rooted map $\mathcal M_\bullet$, and let $\underline{w}$ be the associated double occurrence word. 
	Then,
	\begin{itemize}
	 \item $S$ is the Early DFS of $\mathcal M_\bullet$ if  every  $e\in S$ is a source of  $\vec{\Lambda}(w)$ (or, equivalently, if $w$ does not contain the pattern $f\,e\,f\,e$ with $e\in S$ and $f\notin S$);
	 \item $S$ is the Late DFS of $\mathcal M_\bullet$ 
	 if  every  $e\in S$ is a source of  $\vec{\Lambda}(w)$ (or, equivalently,
	 if $w$ does not contain the pattern $e\,f\,e\,f$ with $e\in S$ and $f\notin S$).
 \end{itemize}
\end{lemma}

\begin{figure}[ht]
	\centering
	\includegraphics[width=.25\textwidth]{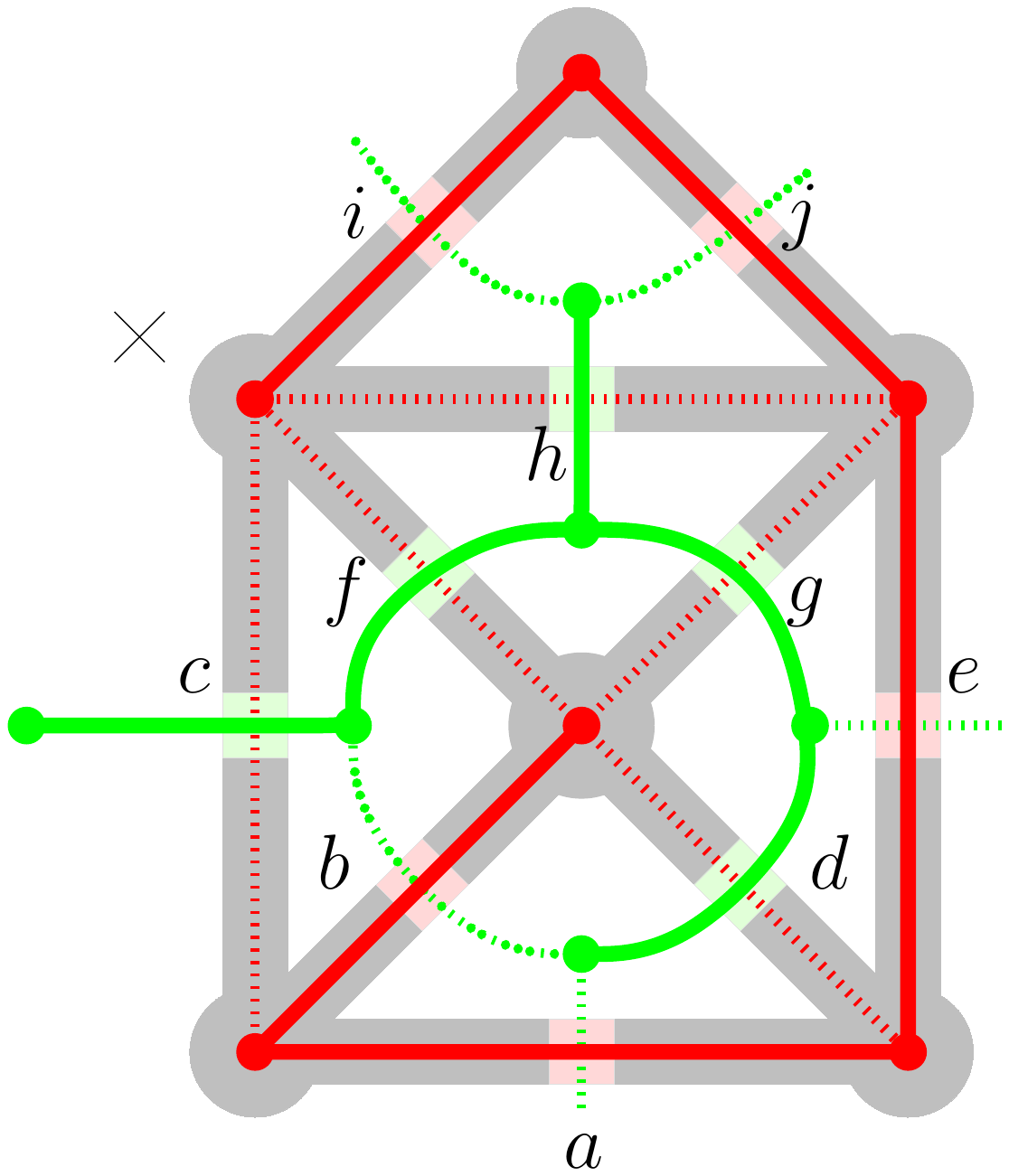}
	\caption{Illustration of the duality between the Early DFS of a planar map and the Late DFS of its dual. Note that because of our definition of the duality of maps, the Early DFS of the primal map $\mathcal M$ corresponds to a DFS giving precedence to the leftmost admissible edge. Re-rooting the dual DFS-tree at $\sigma^{-1}(b_\bullet)$ we get also a DFS giving precedence to the leftmost admissible edge.}
	\label{fig:DFSdual}
\end{figure}

We now deduce a duality results about DFS-trees of planar maps (see \Cref{fig:DFSdual}).

\begin{corollary}
	Let $S$ be a DFS-tree of a rooted $2$-connected planar map $\mathcal M_\bullet$. Then the complement $\overline{S}$ of $S$ is a DFS-tree of $\mathcal M^\ast_\bullet$ (the dual map rooted at the same flag $b_\bullet$) if and only if 
	\begin{itemize}
\item either $S$ is the Early DFS of $\mathcal M_\bullet$ and $\overline{S}$ is the Late DFS of $\mathcal M^\ast_\bullet$,
\item or $S$ is the Late DFS of $\mathcal M_\bullet$ and $\overline{S}$ is the Early DFS of $\mathcal M^\ast_\bullet$.
	\end{itemize}
\end{corollary}
\begin{proof}
 As $\mathcal M$ is planar, $\overline{S}$ is a spanning tree of $\mathcal M^\ast$.

Assume $S$ is the Early (resp. the Late) DFS of $\mathcal M_\bullet$. As the double occurrence word associated to the spanning tree $S$ of $\mathcal M_\bullet$ is the same as the double occurrence word associated to the quasi-tree $\overline{S}$ of $\mathcal M^\ast_\bullet$, we deduce from \Cref{lem:extDFS} that  $\overline S$ is the Late (resp. the Early) DFS of $\mathcal M_\bullet^\ast$. 

Assume 	$\overline{S}$  is a DFS-tree of $\mathcal M^\ast_\bullet$. According to \Cref{lem:DFS}, for every $f\notin S$, either all the $S$-neighbors of $f$ in $\vec{\Lambda}(w)$ are in-neighbors of $f$ (we say that $f$ has type $1$) or  all the $S$-neighbors of $f$ in $\vec{\Lambda}(w)$ are out-neighbors of $f$ (we say that $f$ has type $2$). Dually, for every $e\in S$, either all the $\overline{S}$-neighbors of $e$ in $\vec{\Lambda}(w)$ are in-neighbors of $e$ (we say that $e$ has type $2$) or  all the $\overline{S}$-neighbors of $f$ in $\vec{\Lambda}(w)$ are out-neighbors of $f$ (we say that $e$ has type $1$). Let $\vec \Gamma$ be the bipartite sub-digraph of $\vec{\Lambda}(w)$ induced by the parts $S$ and $\overline{S}$. It is immediate that the neighbors of an edge of type $i$ in $\vec{\Gamma}$ have also the type $i$. As $\mathcal M$ is $2$-connected, the graph $\vec{\Gamma}$ is connected (as easily follows from \Cref{fact:YZ}). Hence, all the edges have the same type. According to \Cref{lem:extDFS}, we deduce that either $S$ is the Early DFS of $\mathcal M_\bullet$ (and $\overline{S}$ is the Late DFS of $\mathcal M^\ast_\bullet$) or $S$ is the Late DFS of $\mathcal M_\bullet$ (and $\overline{S}$ is the Early DFS of $\mathcal M^\ast_\bullet$).
\end{proof}


\section{The Quasi-tree poset}

One can define a poset on the set of all the quasi-trees of a rooted map as follows:
to each pair $(w,I)$, where $w$ is a double occurrence word and $I$ is a subset of letters of $w$ we associate the binary word $b(w,I)$ obtained by extracting from $w$ the subword formed by all the first occurrences of the letters, replacing each letter not in $I$ by $0$ and all letters in $I$ by $1$. Denote by $<_{\rm bin}$ the lexicographic order on the binary strings. Then, a quasi-tree $S$ of $\mathcal M$ is smaller than another quasi-tree $S'$ of $\mathcal M$ is there exists a sequence $S_0=S,S_1,\dots,S_k=S'$ of quasi-trees of $\mathcal M$, where $|S_i\mathbin{\bigtriangleup}S_{i-1}|=2$ 
and $b(w(S_{i-1}), S_{i-1})<_{\rm bin} b(w(S_i),S_i)$
(for $1\leq i\leq k$), where $w(S)$ denotes the first-occurrence word defined by the quasi-tree $S$.
Note that one interpret this poset as a partial order on a pivoting class of (rooted) circle graphs.

\begin{figure}[ht]
	\centering
	\includegraphics[width=.9\textwidth]{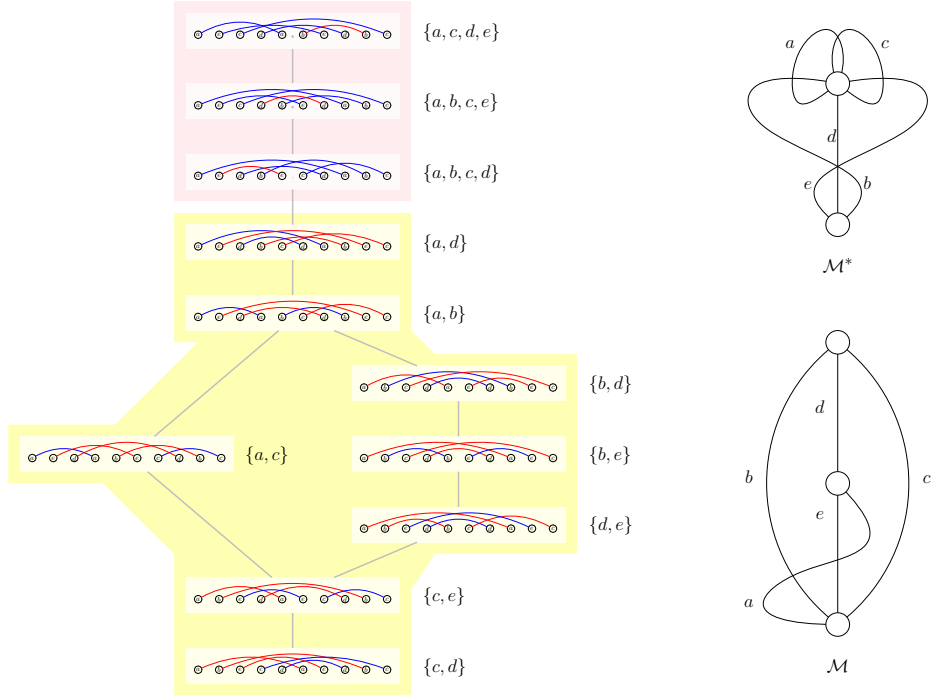}
	\caption{The poset of the quasi-trees of a rooted map. The yellow zone corresponds to spanning trees, while the pink zone corresponds to quasi-trees with genus $1$.
	The minimum element of the poset is the Late DFS-tree, while the maximum one is the complement of the Late DFS-tree of the dual rooted map.}
	\label{fig:nposet}
\end{figure}

As the Late DFS-tree is the minimum of this poset, it is a good candidate to serve as a standard quasi-tree. An immediate property of this tree $T$ is that no tree edge is interlaced on the right. In other words, tree edges are sinks of $\vec{\Lambda}(w)$. On the other hand, if an edge is not in the tree, it is interlaced with at least one tree edge, as we assumed that the map is loopless. Thus, the set of sinks of $\vec{\Lambda}(w)$ dominates  $\vec{\Lambda}(w)$. 

\begin{lemma}
	\label{lem:loopless}
	Double occurrence words $w$ with $|w|=2m$ such that the set of sinks of $\vec{\Lambda}(w)$ dominates  $\vec{\Lambda}(w)$ are in bijection with rooted loopless maps with $m$ edges.
\end{lemma}

\begin{figure}[ht]
	\includegraphics[width=\textwidth]{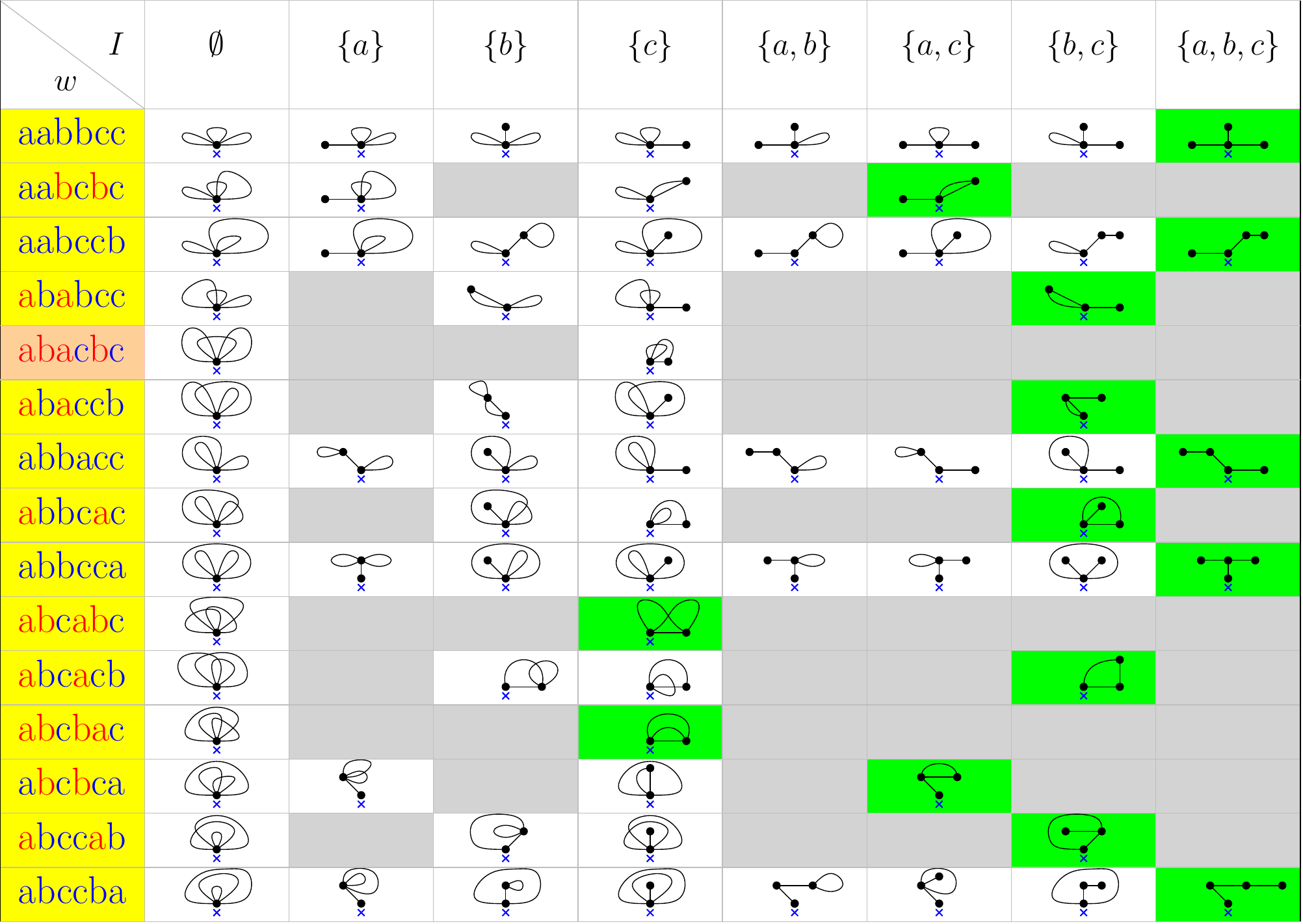}
	\caption{Maps defined by double occurrence words on $3$ symbols. In  $\vec\Lambda(abacbc)$, the only sink is $c$ and it does not dominate $a$. }
	\label{fig:3symb}
\end{figure}

\begin{example}
	There are $14$ orientable rooted loopless maps with $3$ edges, which correspond to the double occurrence sequences distinct from $abacbc$ (see \Cref{fig:3symb}).
\end{example}

\Cref{lem:loopless} can easily be specialized to planar maps.

\begin{lemma}
	\label{lem:looplessP}
	Double occurrence words $w$ with $|w|=2m$ every vertex of $\vec{\Lambda}(w)$ is a source or a sink (i.e. such that no symbol of $w$ is interlaced both on the left and on the right) are in bijection with rooted loopless planar maps with $m$ edges.
\end{lemma}
\begin{proof}
	Consider the bijection between double occurrence words $w$ with $|w|=2m$ such that the set of sinks of $\vec{\Lambda}(w)$ dominates  $\vec{\Lambda}(w)$ and rooted loopless maps with $m$ edges (see \Cref{lem:loopless}).
	In this bijection, the set of all vertices $T(w)$ of  $\vec{\Lambda}(w)$ that have no out-neighbors forms a spanning tree of the map. Thus, the map associated to a word $w$ is planar if and only if the partition $(T(w),V(\vec{\Lambda}(w))\setminus T(w))$ is a partition into two independent sets, that is if and only if no two symbols of $w$ with an out-neighbor in $\vec{\Lambda}(w)$ are adjacent in $\vec{\Lambda}(w)$. 
	This is clearly equivalent to the property that no vertex has both an in-neighbor and an out-neighbor.
\end{proof}

\section{Counting Loopless Rooted Maps}
A \emph{$1$-$2$ occurrence word} is a word in which each letter appears either once or twice. Let $w$ be an $1$-$2$ occurrence word. The symbols that appear only once in $w$ are the \emph{unmatched} symbols of $w$, and they form the set $U(w)$. The symbols that appear twice are the \emph{matched} symbols of $w$.  Removing the unmatched  symbols from $w$, we get the double occurrence word $w^\circ$. An unmatched symbol $a$ is \emph{covered} by a matched symbol $b$ in $w$ if $b\,a\,b$ is a pattern of $w$.

We define two properties for words:
\begin{itemize}
	\item Property \prop{P}: a word $w$ has property \prop{P} if it is a double occurrence word such that the sinks of $\vec\Lambda(w)$ form a dominating set; in other words, every symbol interlaced on the right (non-sink) is interlaced on the right by a symbol that is not interlaced on the right.
	\item Property \prop{Q}: a word $w$ has property \prop{Q} if it is a $1$-$2$ occurrence word, $w^\circ$ has property \prop{P}, and no unmatched symbol of $w$ is covered in $w$ by a sink of  $\vec\Lambda(w)$ (i.e. by a matched symbol of $w$ not interlaced on the right).
\end{itemize}
Note that property \prop{P} obviously implies property \prop{Q} and that 
\cref{lem:loopless} expresses that rooted loopless maps with $n$ edges are in bijection with the double occurrence words with $n$ symbols that satisfy the property \prop{P}.

\begin{lemma}
	\label{lem:decomp}
	Let $w$ be a $1$-$2$ occurrence word.
	
	Then $w$ has property \prop{Q} if and only if 
	\begin{itemize}
		\item either $w$ is empty,
		\item or
		$w=w'\,a$, where $a$ is an unmatched symbol of $w$ and $w'$ has the property \prop{Q}, 
		\item or $w=w_1\,a\,w_2\,a$, where no unmatched symbol of $w$ is in $w_2$, $w_2^\circ$ has property \prop{P}, and $w_1$ has property \prop{Q}.
	\end{itemize} 
	Moreover, such a decomposition, if it exists,  is unique.
\end{lemma}
\begin{proof}
	First, notice that the uniqueness of the decomposition (when it exists) is straightforward.
	We prove the statement by induction on the number of symbols in $w$. If $w$ is empty, the statement is obviously satisfied. If the last symbol $a$ of $w$ is unmatched in $w$ (so $w=w'\,a$) it is immediate that $w$ has property \prop{Q} if and only if $w'$ has property \prop{Q}. So  we can assume that $w$ decomposes as $w=w_1\,a\,w_2\,a$. 
	
	Assume $w$ has property \prop{Q}. As $a$ is  not interlaced on the right (i.e. is a sink of $\vec\Lambda(w)$) we get that no symbol in $w_2$ is unmatched in $w$. That $w_1^\circ$ and $w_2^\circ$ have property \prop{P} follows directly from the hypothesis that $w^\circ$ has property \prop{Q}. Finally, consider an unmatched symbol $f$ of $w_1$ and assume for contradiction that $f$ is covered by a matched symbol $e$ of $w_1$ that is not interlaced on the right in $w_1$ (i.e. is a sink of $\vec\Lambda(w_1)$). Assume first that $e$ is interlaced on the right in $w$ by some matched symbol $g$ of $w$ (i.e. that $e$ is not a sink of $\vec\Lambda(w)$). Then $g$ is interlaced on the right by $a$ in $w$. As this holds for all possible choices of $g$ it follows that $e$ is not dominated by a sink in $\vec\Lambda(w)$, contradicting the assumption that $w$ has property \prop{Q}. Hence, $e$ is a sink of $\vec\Lambda(w)$. As $w$ has property \prop{Q} we get that $f$ is matched in $w$. But then $e$ is interlaced on the right by $f$ in $w$, contradicting the fact that $e$ is a sink of $\vec\Lambda(w)$.  It follows that $w_1$ has property \prop{Q}.
	
	Conversely, assume that $w=w_1\,a\,w_2\,a$,  no unmatched symbol of $w$ is in $w_2$, $w_2^\circ$ has property \prop{P}, and $w_1$ has property \prop{Q}. We first prove the next claim, which determines the unmatched symbols of $w$. 
	\begin{claim}
		A symbol $e$ is a sink of $\vec\Lambda(w)$ if and only if it is a sink of $\vec\Lambda(w_1)$ or  $\vec\Lambda(w_2)$.
	\end{claim}
	\begin{claimproof}
		The left to right implication is straightforward, so we are left with  the right to left implication. 
		
		The symbol $a$ is not interlaced on the right, thus is a sink of $\vec\Lambda(w)$. Every source of $\vec\Lambda(w_2)$ is not interlaced on the right in $w_2$ so cannot be interlaced on the right in $w$ (as the first occurrences of matched symbols of $w$ that are unmatched in $w_2$ are on the left of $w_2$) hence are sinks of $\vec\Lambda(w)$. Consider a sink $e$ of $\vec\Lambda(w_1)$ . If $e$ is not a sink of $\vec\Lambda(w)$, then it covers the first occurrence of a matched symbol $f$ of $w$, whose second occurrence is in $w_2$. But then $e$ covers $f$ in $w_1$ and $f$ is unmatched in $w_1$, contradicting the assumption that $w_1$ has property \prop{Q}.
	\end{claimproof}
	
	Consider a matched symbol $f$ that is not a sink of $w$. If both occurrences of $f$ belong to $w_1$ (resp. to $w_2$) then $f$ is interlaced on the left by a sink of $\vec\Lambda(w_1)$ (resp. a source of $\vec\Lambda(w_2)$), which is (by above claim) a source of $\vec\Lambda(w)$.  Otherwise, if the first occurrence of $f$ belongs to $w_1$ and the second to $w_2$, then $f$ is interlaced on the right by the sink $a$ of $w$. Altogether, we get that $w^\circ$ has property \prop{P}.
	
	Assume for contradiction that some  unmatched symbol of $w$ is covered in $w$ by some sink $e$ of $\vec\Lambda(w)$. According to the assumptions, as $f$ is unmatched in $w$, it belongs to $w_1$. The sink $e$ of $\vec\Lambda(w)$, having its first occurrence in $w_1$ cannot be $a$ or a sink of  $\vec\Lambda(w_2)$. Hence, by the above claim, it is a sink of $\vec\Lambda(w_1)$. It follows that $f$ is covered in $w_1$ by a sink of $\vec\Lambda(w_1)$ contradicting the assumption that $w_1$ has property \prop{Q}.
\end{proof}

Let $G_{n,m}$ denote the number of $1$-$2$ occurrence words with $n$ matched symbols and $m$ unmatched symbols that satisfy the property \prop{Q} if $n,m$ are both non-negative integers, and let $G_{n,m}=0$ otherwise.

\begin{lemma}
	For every  $n\geq 1$ and $m\geq 0$ we have $G_{0,m}=1$ and 
	\begin{equation}
		G_{n,m}=G_{n,m-1}+\sum_{i=0}^{n-1}\sum_{j=0}^{n-1-i}\frac{(2i+j)!}{(2i)!}\,\binom{m+j}{j}\,G_{i,0}\,G_{n-1-i-j,m+j} 
	\end{equation}
\end{lemma}
\begin{proof}
	This is a direct consequence of \cref{lem:decomp}:  every $1$-$2$ occurrence word $w$ with $n$ matched symbol and $m$ unmatched symbols decomposes either as $w'\,a$, where $w'$ as $n$ matched symbols and $m-1$ unmatched symbols, or is obtained (by $w=w_1\,a\,w_2\,a$), for some pair $(i,j)$ of non-negative integers with $i+j\leq n-1$, from a double occurrence word with property \prop{P}  with $i$ symbols ($w_2^\circ$) and an $1$-$2$ occurrence word $w_1$ with $n-1-(i+j)$ matched symbols and $m+j$ unmatched symbols with property \prop{Q} ($w_1$) by inserting $j$ place-holders in $w_2^\circ$ to match $j$ unmatched symbols of $w_1$ (this gives $\binom{2i+j}{j}$ choices), selecting in $w_1$ the $j$ symbols that will be matched in $w_2$ ($\binom{m+j}{j}$ choices) and choosing a matching between these symbols and the place-holders in $w_2$ ($j!$ choices).
\end{proof}

To compute values of $G_{n,0}$ it may be helpful to introduce $M_{n,m}=m!\,G_{n,m}$. (Note that $G_{n,0}=M_{n,0}$.) Then 
$M_{0,m}=m!$ and we have the following recurrence:
\[
M_{n,m}=m\,M_{n,m-1}+\sum_{i=0}^{n-1}\sum_{j=0}^{n-1-i}\binom{2i+j}{j}\,M_{i,0}\,M_{n-1-i-j,m+j}.
\]

Using this recurrence, we computed the values of $G_{n,0}$ for $1\leq n\leq 20$, namely:

$1$,
$3$,
$14$,
$87$,
$672$,
$6204$,
$66719$,
$820395$,
$11370212$,
$175583880$,
$2992513416$,
$55838871492$,
$1132934744671$,
$24846387327825$,
$585953052416226$,
$14791975514747882$,
$398109420366969728$,\linebreak
$11382340640393570304$,
$344600158836813725696$,
$11015256001205535506432$.

Let $s=2n+m, t=n+m$, and $N_{s,t}=M_{s-t,2t-s}$.
Then we have
\[
N_{s,t}=(2t-s)N_{s-1,t-1}+ \sum_{k=0}^{s-t-1}\sum_{\ell=k}^{2k}\binom{\ell}{k}\, N_{2(\ell-k),\ell-k}N_{s-2-\ell, t-1}.
\]
\section{Back to the Planar Case}
Interestingly, a variant of the problem of enumerating  rooted maps in which only separating loops are forbidden was solved by Walsh and Lehman \cite{walsh1975counting} for any given genus. As every loop of a planar map is separating,  it results from their work that the number of loopless rooted planar maps with $n$ edges is given by the closed expression $\frac{2\,(4n+1)!}{(n+1)!\,(3n+2)!}=\binom{4n+1}{n+1} - 9\binom{4n+1}{n-1}$ (sequence A000260 of the On-line encyclopedia of integer sequences). We refer the interested reader to  \cite{BENDER1985235} for more results on enumerations of loopless planar maps.

In this section, we consider how this problem is connected to the bijections we defined in the general case.
We define two new properties for words:
\begin{itemize}
	\item Property \prop{N}: a word $w$ has property \prop{N} if it is a double occurrence word such that every vertex of $\vec\Lambda(w)$ is a source or a sink or, equivalently, if no symbol of $w$ is interlaced both on the left and on the right;
	\item Property \prop{N'}: a word $w$ has property \prop{N'} if it is a $1$-$2$ occurrence word, $w^\circ$ has property \prop{N}, and every symbol covering an unmatched symbol of $w$ is a sink of $\vec\Lambda(w)$ (i.e. is not interlaced on the right).
\end{itemize}
Note that property \prop{N} obviously implies property \prop{N'} and that 
\cref{lem:looplessP} expresses that rooted loopless planar maps with $n$ edges are in bijection with the double occurrence words with $n$ symbols that satisfy the property \prop{N}.

For a word $w$, we denote by $\widetilde{w}$ the \emph{reversal} of $w$, that is the word obtained from $w$ by reversing the order of the letters.

\begin{figure}[ht]
		\centering
	\includegraphics[width=.9\textwidth]{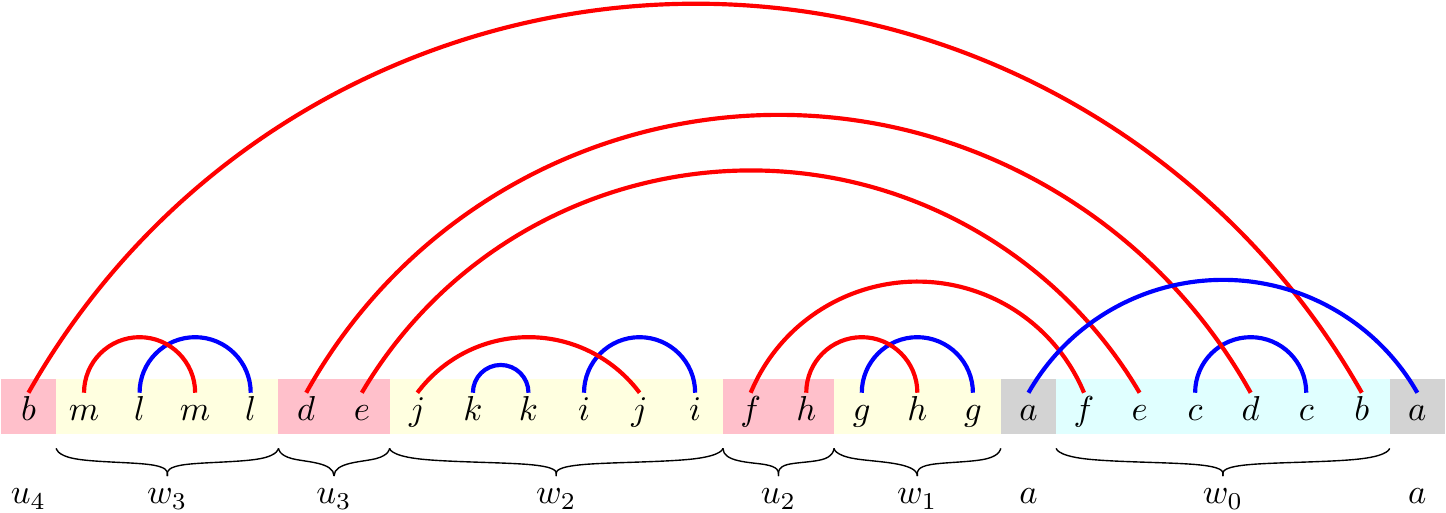}
	\caption{Example of a decomposition of a double occurrence word with property \prop{N}. The words $u_1$ and $w_4$ are empty.}
	\label{fig:decompN}
\end{figure}

\begin{lemma}
	Let $w$ be a non-empty $1$-$2$ occurrence word.
	Then $w$ has property $\prop{N'}$ if and only if 
	\begin{itemize}
		\item 	either $w$ is empty,
		\item or $w=w'\,a$ where $a$ is an unmatched symbol of $w$ and $w'$ has property \prop{N'}, 
		\item or $w=w_k\,u_k\dots w_1\,u_1\,a\,w_0\,a$, where `
		\begin{itemize}
		\item $a$ is a symbol,
		\item the words $u_1,\dots,u_k$ contain (together) exactly the first occurrences of the matched symbols of $w$ that belong to $w_0$ but are unmatched in $w_0$,
		\item the words $u_2,\dots,u_k$  and $w_1,\dots,w_{k-1}$ are not empty,
		\item  the pattern $z$ of $w_0$ with all symbols in $U(w_0)\setminus U(w)$ satisfy $u_k\,\dots\,u_1=\widetilde{z}$, 
		\item every symbol in $z$ occurs before all the elements of $U(w_0)\cap U(w)$,
		\item  $w_0$ and $w_k$ have property \prop{N'},
		\item every $w_i$ with $1\leq i<k$ has property \prop{N}.		 
	 \end{itemize}
 
	\end{itemize}	
	
	Moreover, the decomposition, if it exists, is unique.
\end{lemma}
\begin{proof}
	Let $w$ be a non-empty $1$-$2$ occurrence word whose last symbol is matched in $w$ and let 	$w=w_k\,u_k\dots w_1\,u_1\,a\,w_0\,a$, where $a$ is a symbol, $u_1,\dots,u_k$ contain the first occurrences of the matched symbols of $w$ that belong to $w_0$ but are unmatched in $w_0$, $u_2,\dots,u_k$ are not empty, and $w_i$ is not empty if $1\leq i<k$.  Note that this decomposition is uniquely defined. Let $z$ be the  pattern of $w_0$ with all symbols in $U(w_0)\setminus U(w)$.
	
	Assume $w$ has property \prop{N'}. 
	As every symbol in $z$ is interlaced on the right by $a$, no symbol in $z$ is interlaced on the left, thus no symbol in $u_1,\dots,u_k$ is covered in
	$w_k\,u_k\dots w_1\,u_1$.
As $u_k$ is not empty, there exists a symbol $f$ with second occurrence in $w_0$ and first occurrence in $u_k$. The symbol $f$ is interlaced on the right by $a$ in $w$. By property \prop{N'}, $f$ does not cover any unmatched symbol. In particular, $w_{k-1},\dots,w_{1}$ do not contain any unmatched symbol of $w$. Moreover no symbol $f$ can have a first occurrence in $w_j$ and a second occurrence in $w_i$ for $1\leq i<j\leq k$ as otherwise it would interlace on the left a symbol $e$ with second occurrence in $w_0$ and first occurrence in $u_j$ (such a symbol exists as 
	$u_j$ is not empty), which is itself interlaced on the right by $a$, contradicting the hypothesis that $w^\circ$ has property \prop{N}. It follows that all the $w_i$ with $1\leq i<k$ are double occurrence words.
	It is straightforward that if $w'$ is a pattern of $w$ with $U(w')\subseteq U(w)$ then $w'$ has property \prop{N'}. Hence we get  $u_k\dots u_0=\widetilde{z}$ (by considering the pattern including the two occurrences of $a$ of the symbols in $z$), that $w_i$ has property \prop{N} for $1\leq i<k$, and that $w_k$ has property \prop{N'}. Also, assume that some $e\in U(w)\cap U(w_0)$ appears before some $f\in U(w_0)\setminus U(w)$ in $w_0$. Then $f$ covers $e$ in $w$ and $f$ is interlaced on the right by $a$ in $w$, contradicting the hypothesis that $w$ has property \prop{N'}.
	We now prove that $w_0$ has property \prop{N'}.
	First note that  every symbol in $U(w)\cap U(w_0)$ is not covered in $w_0$ by a symbol that is interlaced on the  right in $w_0$, for otherwise the same would hold in $w$, contradicting property \prop{N'}. 
	Let $f\in U(w_0)\setminus U(w)$ and assume for contradiction that $f$ is covered by $e$ in $w_0$ that is interlaced on the right by $g$. Then $e$ is interlaced (in $w$) on the right by $g$ and on the left by $f$, contradicting the hypothesis that $w$ has property \prop{N'}.

	Conversely, assume that $u_k\dots u_1=\widetilde{z}$ and that $w_0,w_k$ have property \prop{N'}, that all $w_i$ with $1\leq i<k$ have property \prop{N}, and that all the elements in $U(w_0)\cap U(w)$ are on the left of all the elements in $U(w_0)\setminus U(w)$.
	Assume for contradiction that $w$ does not have the property \prop{N'}. Then either there exists a matched symbol $f$ interlaced both on the left and on the right, or there exists an unmatched symbol $f$ covered by a symbol $e$ that is interlaced on the right. 
	Assume that there exists a matched symbol $f$ interlaced both on the left and on the right. Obviously $f\neq a$.
	Also, $f$ is not a symbol in $z$ as these symbols are interlaced only on the right. If $f$ belongs to $w_i$ with  $1\leq i\leq k$ then the two occurrences of $f$ as well as the two occurrences of every symbol interlaced with $f$ in $w$ belong to $w_i$, contradicting the assumption that $w_i$ has property \prop{N'}. Hence the two occurrences of $f$ belong to $w_0$. As $w_0$ has property \prop{N'} the symbol $f$ is not interlaced both on the left and on the right in $w_0$. It follows that $f$ is interlaced on the right in $w_0$ and is interlaced on the left by some symbol $g$ in $z$, contradicting property \prop{N'}. 
	Now assume that there exists an unmatched symbol $f$ covered by a symbol $e$ that is interlaced on the right. Then either $f$ belongs to $w_0$ or $f$ belongs to $w_k$. If $f$ belongs to $w_k$ so do the two occurrences of $e$ as well as the two occurrences of any symbol interlaced with $e$, contradicting the hypothesis that $w_k$ has property \prop{N'}. If $f$ belongs to $w_0$ then either the two occurrences of $e$ (as well as the two occurrences of the symbol interlaced with $e$ on the left) belong to $w_0$, contradicting the assumption that $w_0$ has property \prop{N'}, or $e\in U(w_0)\setminus U(w)$, contradicting the hypothesis that all the elements of $U(w_0)\cap U(w)$ (like $f$) are on the right of all the elements in $U(w_0)\setminus U(w)$ (like the second occurrence of $e$).
\end{proof}

Remark that the condition  $u_k\dots u_1=\widetilde{z}$ and the condition that all the elements in $U(w_0)\cap U(w)$ are on the right of all the elements in $U(w_0)\setminus U(w)$ simply mean that we do not have any choice when we will match unmatched symbols of $w_0$, what will simplify the counting.

Let $T_{n,m}$ be the number of $1$-$2$ occurrence words with  $n$ matched symbols and $m$ unmatched symbols that satisfy property \prop{N'}. 
\begin{lemma}
	Let
	\[
	F(x,y)=\sum_{n\geq 0}\sum_{m\geq 0}T_{n,m}\,x^n\,y^m.
	\]
	Then $F(x,y)$ is a solution of the equation
	\begin{equation}
		F(xy,y)=1+yF(xy,y)+xy\frac{F(xy,y)^2}{1-yF(xy,0)}
	\end{equation}
\end{lemma}
\begin{proof}
	The first summand is $1$, for the empty word; the second is $yF(xy,y)$ for $w'$ followed by $a$; the third 
	summand corresponds to a matched $a$ ($xy$), the word $w_k$ with property \prop{N'} ($F(xy,y)$), a 	
	sequence ($(1-(F(xy,0-1)y(1-y)^{-1}))^{-1}$) of words formed
	a non-empty sequence of unmatched symbols ($y(1-y)^{-1}$) followed by
	 by a non-empty words with property \prop{N} ($F(xy,0)-1$),  
	a sequence of unmatched symbols ($(1-y)^{-1}$),  and finally the word $w_0$ with property \prop{N'} ($F(xy,y)$).
\end{proof}

\providecommand{\bysame}{\leavevmode\hbox to3em{\hrulefill}\thinspace}
\providecommand{\MR}{\relax\ifhmode\unskip\space\fi MR }
\providecommand{\MRhref}[2]{%
	\href{http://www.ams.org/mathscinet-getitem?mr=#1}{#2}
}
\providecommand{\href}[2]{#2}

\end{document}